\documentstyle[12pt,fleqn]{article}

\def\thebibliography#1{\section{References}\list
  {\arabic{enumi}.}{\settowidth\labelwidth{#1}\leftmargin\labelwidth
    \advance\leftmargin\labelsep
    \usecounter{enumi}}
    \def\newblock{\hskip .11em plus .33em minus .07em}
    \sloppy\clubpenalty4000\widowpenalty4000
    \sfcode`\.=1000\relax}

\newcommand{\bx}[1]{{\rm\fbox{$ #1$}}}

\def\op#1{\mathop{\fam0 #1}\limits}

\newcommand{\pr}{{\rm pr}}
\newcommand{\Id}{{\rm Id}}
\newcommand{\Ker}{{\rm Ker\,}}
\newcommand{\nm}[1]{\mid {#1}\mid}
\newcommand{\df}{{\rm def}}

\newcommand{\ry}{{\rm y}}

\newcommand{\cO}{{\cal O}}
\newcommand{\cT}{{\cal T}}

\newcommand{\cL}{{\cal L}}
\newcommand{\cV}{{\cal V}}
\newcommand{\cE}{{\cal E}}
\newcommand{\cH}{{\cal H}}
\newcommand{\cF}{{\cal F}}
\newcommand{\cC}{{\cal C}}
\newcommand{\bL}{{\bf L}}
\newcommand{\R}{{\bf R}}
\newcommand{\al}{\alpha}
\newcommand{\bt}{\beta}
\newcommand{\dl}{\delta}
\newcommand{\la}{\lambda}

\newcommand{\bla}{{\bf\Lambda}}
\newcommand{\bom}{{\bf\Omega}}
\newcommand{\bth}{{\bf\Theta}}
\newcommand{\f}{\phi}

\newcommand{\vt}{\vartheta}

\newcommand{\p}{\pi}
\newcommand{\bpi}{{\bf\Pi}}

\newcommand{\om}{\omega}
\newcommand{\Om}{\Omega}
\newcommand{\m}{\mu}

\newcommand{\g}{\gamma}
\newcommand{\G}{\Gamma}
\newcommand{\e}{\epsilon}
\newcommand{\ve}{\varepsilon}
\newcommand{\up}{\upsilon}
\newcommand{\th}{\theta}

\newcommand{\si}{\sigma}
\newcommand{\Si}{\Sigma}

\newcommand{\Y}{Y\to X}
\newcommand{\w}{\wedge}
\newcommand{\wt}{\widetilde}
\newcommand{\wh}{\widehat}
\newcommand{\ol}{\overline}
\newcommand{\dr}{\partial}
\newcommand{\fl}{\flat}
\newcommand{\sh}{\sharp}

\newcommand{\lra}{\op\longrightarrow}

\newcommand{\llra}{\longleftrightarrow}
\newcommand{\har}{\op\hookrightarrow}
\newcommand{\xx}{\times}
\newcommand{\ot}{\otimes}

\newcommand{\ap}{\approx}

\let\ssection=\section
\renewcommand{\section}{\setcounter{equation}{0}\ssection}

\newcounter{eqalph}[section]
\newcounter{equationa}[section]
\newcounter{example}[section]
\newcounter{remark}[section]
\newcounter{theorem}[section]
\newcounter{proposition}[section]
\newcounter{lemma}[section]
\newcounter{corollary}[section]
\newcounter{definition}[section]

\def\theremark{\arabic{section}.\arabic{remark}}

\def\thedefinition{\arabic{section}.\arabic{definition}}

\newenvironment{imp}{\bigskip{ \huge $\wr$}}{{ \huge $\wr$}\bigskip} 
\newenvironment{proof}{\noindent {\bf Proof.}}{{\Large $\bullet$} \medskip }
\newenvironment{example}{\refstepcounter{remark} \medskip\noindent{\bf Example
\theremark}.}{{\Large $\bullet$} \medskip }
\newenvironment{remark}{\refstepcounter{remark} \medskip\noindent{\bf Remark
\theremark}.}{{\Large $\bullet$} \medskip }
\newenvironment{theorem}{\refstepcounter{definition} \medskip\noindent{\sc
Theorem \thedefinition}.}{$\Box$\medskip }
\newenvironment{proposition}{\refstepcounter{definition} \medskip\noindent{\sc
Proposition \thedefinition}.}{$\Box$ \medskip }
\newenvironment{lemma}{\refstepcounter{definition} \medskip\noindent{\sc Lemma
\thedefinition}.}{ $\Box$\medskip }
\newenvironment{corollary}{\refstepcounter{definition} \medskip\noindent{\sc
Corollary \thedefinition}.}{ $\Box$\medskip }
\newenvironment{definition}{\refstepcounter{definition} \medskip\noindent{\sc
Definition \thedefinition}.}{$\Box$ \medskip }

\newcommand{\beq}{\begin{equation}}
\newcommand{\eeq}{\end{equation}}
\newcommand{\ben}{\begin{eqnarray}}
\newcommand{\een}{\end{eqnarray}}
\newcommand{\be}{\begin{eqnarray*}}
\newcommand{\ee}{\end{eqnarray*}}
\newcommand{\bea}{\begin{eqalph}}
\newcommand{\eea}{\end{eqalph}}

\newenvironment{eqalph}{\stepcounter{equation}
\setcounter{equationa}{\value{equation}}
\setcounter{equation}{0}

\begin{eqnarray}}{\end{eqnarray}
\setcounter{equation}{\value{equationa}}}

\begin{document}

\hbox{}

\pagestyle{empty}

\vskip 1cm

\begin{center}

{\huge Differential Geometry 
\bigskip

of Time-Dependent Mechanics}
\vskip 1cm

{\Large \bf G. Giachetta, L. Mangiarotti,}
\medskip

Department of Mathematics and Physics, University of Camerino,

62032 Camerino (MC), Italy 

E-mail: mangiaro@camserv.unicam.it
\bigskip

{\Large \bf G. Sardanashvily}
\medskip

Department of Theoretical Physics, Moscow State University

117234 Moscow, Russia

E-mail: sard@grav.phys.msu.su
\end{center}

\vskip 1cm

\begin{abstract}
The usual formulations of time-dependent mechanics start from a given 
splitting $Y=\R\times M$ of the coordinate bundle $Y\to\R$. From physical
viewpoint, this splitting means that a reference frame has been chosen. 
Obviously, such a splitting is broken under reference frame transformations 
and time-dependent canonical transformations. Our goal is to formulate
time-dependent mechanics in gauge-invariant form, i.e., independently of any
reference frame. The main ingredient in this formulation is a connection on
the bundle $Y\to\R$ which describes an arbitrary reference frame. We 
emphasize the following peculiarities of this approach to time-dependent
mechanics. A phase space does not admit any canonical contact or 
presymplectic structure which would be preserved under reference frame
transformations, whereas the canonical Poisson  structure is degenerate. A
Hamiltonian fails to be a function on a phase space. In particular, it can 
not participate in a Poisson bracket so that the evolution equation is not
reduced to the Poisson bracket. This fact becomes relevant to the quantization
procedure. Hamiltonian and Lagrangian formulations  of time-dependent
mechanics are not equivalent. A degenerate Lagrangian admits a set of
associated Hamiltonians, none of which describes the whole mechanical system
given by this Lagrangian. 
\end{abstract}

\newpage

\tableofcontents 

\newpage

\pagestyle{headings}
\markright{}{}

\section{Introduction}

There is an extensive literature both on autonomous
\cite{abra,arno,guil,libe} and time-dependent mechanics
\cite{binz,cari89,eche,leon,mass,mora} (this list of references is of course
far from being exhaustive). The mechanics of autonomous systems
is phrased in terms of symplectic geometries on even-dimensional manifolds,
in particular, on the cotangent bundle $T^*M$ of a manifold $M$. At the same
time, the usual formulations of time-dependent mechanics are developed on 
$\R\xx T^*M$. From physical viewpoint, this means that some
reference frame has been choosen.

In this paper, our goal is the formulation of time-dependent mechanics in 
gauge-invariant form, i.e., independently of any reference frame. The main
ingredient in this formulation is a bundle $Y\to X$ over a 1-dimensional
base $X$. In such a context, a complete connection on $\Y$ defines a
reference frame. From the mathematical viewpoint, a complete connection on
$\Y$ is equivalent to give a splitting $Y\simeq X\xx M$.
This, in turn, implies the splitting of the covariant phase space
$\Pi=V^*Y\simeq X\xx T^*M$ (where 
$V^*Y$ denotes the vertical cotangent bundle of $\Y$).

We emphasize the following  peculiarities of time-dependent mechanics.

(i) The phase space does
not admit any canonical contact or presymplectic structure which
would be maintained under changes of reference frames. At the same time, we
have the canonical Poisson structure, but it is necessarily degenerate. 

(ii) The Hamiltonian fails to be a function on a
 phase space (see (\ref{m002})). In particular, it can not participate in a
Poisson bracket. It follows that the evolution equation is not reduced to the
Poisson bracket. As a consequence, integrals of motion can not be defined as
functions in involution with the Hamiltonian. 

(iii) Canonical transformations fail to admit
even local generating functions in general.

(iv) The spray evolution equation is not maintained under general
reference frame transformations.

(v) Hamiltonian and Lagrangian
formulations of time-dependent mechanics are not equivalent. 
A degenerate Lagrangian admits a set of 
associated Hamiltonians none of which describes the whole mechanical system
given by this Lagrangian. At the same time, we have not any canonical Poisson
or contact structure on a configuration space of Lagrangian mechanics.

We develop time-dependent mechanics as the particular
case of field theory on bundles $Y\to X$ over a
$n$-dimensional base \cite{binz,cari91,sard95}. In this aproach, the physical
variables are described by sections of $\Y$ (where
$\dim X>1$ in field theory and  $\dim X=1$ in mechanics).
If $n>1$, the Hamiltonian partner of the
first order Lagrangian machinery is the polysymplectic Hamiltonian formalism
\cite{cari,giac95,gunt,kijo,sard94,sard95}. If $n=1$, we show that 
this formalism provides the differential geometric description of
time-dependent Hamiltonian mechanics. In particular, the $n=1$ polysymplectic
Hamiltonian form is exactly the integral invariant of Poincar\'e--Cartan 
\beq
H=p_idy^i -\cH dt \label{m001}
\eeq
\cite{arno,libe} on the  phase space $V^*Y$ coordinatized by $(t,y^i,p_i)$.
The form (\ref{m001}) is written in gauge-invariant form
where the Hamiltonian $\cH$ is not a function, but a section of the affine
bundle $T^*Y\to V^*Y$.  The bundle $T^*Y\to X$ coordinatized by
$(t,y^i,p,p_i)$ plays the role of the  phase space in the homogeneous
formalism.  We have the splitting of a Hamiltonian
\beq
\cH=p_i\G^i +\wt\cH, \label{m002}
\eeq
where $\wt\cH$ is a Hamiltonian function and $\G$ is a connection on
$Y\to X$ describing a (local) reference frame.

By analogy with field theory, we
may talk of {\it sui generis} gauge-invariant formulation of mechanics.
Such formulation enables us both to describe a mechanical system without a
preferable reference frame (e.g., in relativistic mechanics and gravitation
theory) and  to analyze phenomena which depend assentially on the choice of a
reference frame (e.g., the energy-momentum conservation laws). Also
quantizations with respect to different reference frames are not equivalent.

\newpage

\section{Preliminaries}

This Section includes the main notions of differential geometry and jet
formalism which we need in sequel. 

{} From a pragmatic viewpoint, we widely use coordinate expressions, but all
objects satisfy the corresponding transformation 
laws and are globally defined.

Throughout, morphisms  are smooth mappings of
class $C^\infty$. Manifolds are real, finite-dimensional, 
paracompact and connected. 

\begin{remark}\label{onedim} The only 1-dimensional manifolds
obeying these conditions are the real line $\R$ and the circle $S^1$.
\end{remark}

We use the symbols $\ot$, $\vee$ and $\w$ for tensor, symmetric and
exterior products respectively. By $\rfloor$ is meant the contraction of
multivectors and differential forms. 
The natural projections of the product $A\times B$ are denoted by
\be
\pr_1:A\times B\to A, \qquad  \pr_2:A\times B\to B.
\ee

Let $Z$ be a manifold coordinatized by $(z^\la)$. The tangent bundle $TZ$ and
the cotangent bundle $T^*Z$ of $Z$  are equipped with the induced coordinates
$(z^\la, \dot z^\la)$ and $(z^\la, \dot z_\la)$
relative to the holonomic fibre bases $\{\dr_\la\}$ and
$\{dz^\la\}$ for $TZ$ and $T^*Z$ respectively.
By $Tf: TZ\to TZ'$ is meant the tangent morphism to morphism $f:Z\to Z'$.

We recall here the following kinds of morphisms:  
immersion, imbedding, submersion, and projection. A morphism
$f:Z\to Z'$ is called {\it immersion} if the tangent morphism $Tf$ at every
point
$z\in Z$ is an injection. When $f$ is both an immersion and an injection, its
image is said to be a submanifold of $Z'$. A submanifold which also is a
topological subspace is called imbedded submanifold.
A mapping $f:Z\to Z'$ is called {\it submersion} if the tangent morphism $Tf$ 
at every point $z\in Z$
is a surjection. If $f$ is both a submersion and a surjection, it is
termed {\it projection} or {\it fibred manifold}.

\subsection{Bundles}

Let $\pi :\Y$ be a fibred manifold over the base $X$.
It is provided with an atlas of fibred coordinates $(x^\la, y^i)$, where
$(x^\la)$ are coordinates of $X$. 

Hereafter, by a {\it bundle} is meant a locally trivial fibred manifold, i.e.
there exists an open covering $\{U_\xi\}$ of $X$ and local diffeomerphisms
\be
\psi_\xi :\pi^{-1}(U_\xi)\to U_\xi\times V,
\ee
where $V$ is the standard fibre of $Y$. The collection
\be
\Psi = \{U_\xi, \psi_\xi, \rho_{\xi\zeta}\}, \qquad
\psi_\xi (y)=(\rho_{\xi\zeta}\circ\psi_\zeta)(y),
\qquad y\in \pi^{-1}(U_\xi\cap U_\zeta),
\ee
of the splittings $\psi_\xi$ together with the transition functions
$\rho_{\xi\zeta}$ constitute a {\it bundle atlas} of $Y$. The associated 
bundle coordinates of $Y$ are
\be
y^i(y)=(v^i\circ\pr_2\circ \psi_\xi)(y), \qquad \pi (y)\in U_\xi, 
\ee
where $(v^i)$ are fixed coordinates of the standard fibre $V$ of $Y$. 
A bundle $\Y$ is said to be {\it trivlaizable} 
if it admits a global splitting
$Y\simeq X\xx V$. Different such splittings differ from each other in 
projections of $Y$ onto $V$.

\begin{theorem}\label{theo1}
Every bundle over a contractible paracompact manifold is trivializable.
\end{theorem}

\begin{theorem}\label{theo2} If the standard fibre of a bundle over 
a paracompact base is diffeomorphic to $\R^m$, this bundle 
has a global section.
\end{theorem}

By a {\it bundle morphism} of $Y\to X$ to $Y'\to X'$ is meant a pair of
mappings $(\Phi,f)$ which form the commutative diagram
\be
\begin{array}{rcccl}
& {Y} &  \op\longrightarrow^{\Phi} & {Y'} &  \\
{} &\put(0,10){\vector(0,-1){20}} & & \put(0,10){\vector(0,-1){20}} & {} \\
& {X} & \op\longrightarrow_{f} & {X'} &
\end{array}
\ee
One says that $\Phi$ is a bundle morphism over $f$ (or over $X$ if
$f=\Id_X$). 

Given a bundle $Y\to X$, every mapping
$f: X' \to X$ yields a bundle  $f^*Y$ over $X'$ which is called
{\it pullback} of the bundle $Y$ by $f$. The
fibre of $f^*Y$ at a point $x'\in X'$ is that of $Y$ at the
point $f(x)\in X$. In particular, the product of bundles
$\pi:\Y$ and $\pi':Y'\to X$ over $X$ is the  pullback
\be
\pi^*Y'={\pi'}^*Y=Y\op\xx_X Y'.
\ee
For the sake of simplicity, we shall denote the pullbacks
$Y\op\times_X TX$, $Y\op\times_X T^*X$
of tangent and cotangent bundles of $X$ by $TX$ and $T^*X$.

\begin{remark}\label{rem1}
Let $\pi:Y\to X$ be a bundle. Every diffeomorphism $\rho$ of a manifold
$Y$ which does not preserve the fibration $\pi$ defines a new bundle
$\pi\circ\rho^{-1}:Y\to X$. Obviously, $\rho$ is 
isomorphism of the bundle $\pi$ to the bundle $\pi\circ\rho^{-1}$ over $X$.
At the same time, fibrations $\pi$ and $\pi\circ\rho^{-1}$ of $Y$ are not
equivalent. Let 
$\rho$ be a bundle isomorphism of $\Y$ over $X$. Given an atlas $\Psi =
\{U_\xi, \psi_\xi\}$ of $Y$, there always exists the atlas 
\beq
\Psi\circ \rho^{-1} = \{U_\xi,
\psi'_\xi=\psi_\xi\circ\rho^{-1}\} \label{m136}
\eeq
of $Y$ such that the bundle coordinates of points
$\rho(y)$ with respect to $\Psi\circ \rho^{-1}$ coincide with the bundle
coordinates  of points $y$ with respect to $\Psi$.
\end{remark}

\subsection{Differential forms and multivectors}

In this Section $Z$ is an $m$-dimensional manifold.

An exterior differential {\it $r$-form} (or simply a $r$-form) $\f$ on a
manifold $Z$ is a section of the bundle $\op\w^r T^*Z\to Z$. The 1-forms
are called the {\it Pfaffian forms}. We utilize the coordinate expression
\be
\f =\f_{\la_1\dots\la_r}
dz^{\la_1}\w\cdots\w dz^{\la_r}, \qquad \nm\f=r,
\ee
where the summation is over all ordered collections $(\la_1,...,\la_r)$.  
We denote by $\cO^r(Z)$ and $\cO(Z)$ the vector space of
$r$-forms and the $\bf Z$-graded algebra of all differential forms on a
manifold $Z$ respectively. 

Given a map $f:Z\to Z'$, by $f^*\f$ is meant the pullback on $Z$ 
of a form $\f$ on $Z'$ by $f$. We recall the relations
\be
f^*(\f\w\si) =f^*\f\w f^*\si, \qquad  df^*\f =f^*(d\f). 
\ee

{\it Contraction} of a vector field $u = u^\m\dr_\m$ and a r-form $\f$
on $Z$ is given in coordinates by 
\be
u\rfloor\f = u^\m
\f_{\m\la_1\ldots\la_{r-1}} 
dz^{\la_1}\w\cdots \w dz^{\la_{r-1}}.
\ee
There is the relation
\be
u\rfloor(\f\w\si)= u\rfloor\f\w\si +(-1)^{\nm\f}\f\w u\rfloor\si.
\ee

The {\it Lie derivative} ${\bf L}_u\f$ of an exterior form $\f$ along
a vector field $u$ is defined to be
\be
\bL_u\f =u\rfloor d\f +d(u\rfloor\f).
\ee
It satisfies the relation
\be
\bL_u(\f\w\si)= \bL_u\f\w\si +\f\w\bL_u\si.
\ee

\begin{example}\label{sympl} Let $\Om$ be a 2-form
on $Z$. It defines the bundle morphism 
\beq
\bx{\Om^\fl: TZ\to T^*Z, \qquad \Om^\fl(v)\op=^\df- v\rfloor\Om(z),
\qquad v\in T_zZ.} \label{m52}
\eeq
In coordinates, if $\Om=\Om_{\m\nu}dz^\m\w dz^\nu$ and $v=v^\m\dr_\m$, then
\be
\Om^\fl(v)= -\Om_{\m\nu} v^\mu dz^\nu.
\ee
One says that $\Om$ is of constant rank $k$ if the corresponding morphism
(\ref{m52}) is of constant rank $k$ (i.e., $k$ is the greatest integer $n$
such that $\Om^n$ is not the zero form). The rank of a {\it nondegenerate}
2-form is equal to $\dim Z$. A nondegenerate closed 2-form is
called the {\it symplectic form}. 
\end{example}

A {\it multivector field} $\vt$ of degree $r$ (or simply a $r$-vector
field) on a manifold $Z$ is a section of the bundle $\op\w^r TZ\to
Z$. It is given by the coordinate expression
\be
\vt =\vt^{\la_1\dots\la_r}
\dr_{\la_1}\w\cdots\w \dr_{\la_r}, \qquad \nm\vt =r,
\ee
where summation is over all ordered collections $(\la_1,...,\la_r)$.

We denote by $\cV^r(Z)$ and $\cV(Z)$ the vector space of 
$r$-vector fields and the $\bf Z$-graded algebra of all 
multivector fields on a manifold $Z$ respectively. The latter is provided
with the {\it Schouten-Nijenhuis bracket} 
\be
[.,.]: \cV^r(Z)\xx\cV^s(Z) \to \cV^{r+s-1}(Z)
\ee
which generalizes the Lie
bracket of vector fields \cite{bhas,vais}. This bracket has the
coordinate expression
\be
&& \vt=\vt^{\la_1\ldots\la_r}\dr_{\la_1}\w\cdots\w\dr_{\la_r}, \qquad 
\up=\up^{\al_1\ldots\al_s}\dr_{\al_1}\w\cdots\w\dr_{\al_s}, \\
&& [\vt,\up] =\vt\star\up +(-1)^{\nm\vt\nm\up}\up\star\vt, \\
&& \vt\star\up = \vt^{\mu\la_1\ldots\la_{r-1}}
\dr_\m\up^{\al_1\ldots\al_s}\dr_{\la_1}\w
\cdots\w\dr_{\la_{r-1}}\w\dr_{\al_1}\w\cdots\w\dr_{\al_s}.
\ee
There are the relations
\be
&&[\vt,\up]=(-1)^{\nm\vt\nm\up}[\up,\vt], \\
&&[\nu,\vt\w\up] =[\nu,\vt]\w\up +(-1)^{\nm\nu\nm\vt +\nm\vt}\vt\w[\nu,\up],
\\ 
&& (-1)^{\nm\nu\nm\vt +\nm\nu}[\nu,\vt\w\up] +(-1)^{\nm\vt\nm\nu
+\nm\vt}[\vt,\up\w\nu] +(-1)^{\nm\up\nm\vt +\nm\up}[\up,\nu\w\vt] =0.
\ee

\begin{example}\label{biv}
Let $w=w^{\m\nu}\dr_\m\w\dr_\nu$ be a {\it bivector field}. We have 
\beq
[w,w] = w^{\m\la_1}\dr_\m
w^{\la_2\la_3}\dr_{\la_1}\w\dr_{\la_2}\w\dr_{\la_3}. \label{m101}
\eeq
Every bivector field $w$ on a manifold $Z$ yields
the associated bundle morphism $w^\sh: T^*Z\to TZ$ defined by
\beq
\bx{w^\sh(p)\rfloor q \op=^\df w(z)(p,q),
\qquad w^\sh(p)= w^{\m\nu}(z)p_\m\dr_\nu, \qquad p,q\in T^*_zZ.} \label{m51}
\eeq
A bivector field $w$ whose bracket (\ref{m101}) vanishes is
called the {\it Poisson bivector field}.
\end{example}

Elements of the tensor product $\cO^r(Z)\ot\cV^1(Z)$ 
are called the {\it tangent-valued $r$-forms} on $Z$. They are sections
\be
\phi = \phi_{\la_1\ldots\la_r}^\m dz^{\la_1}\w\cdots\w dz^{\la_r}\ot\dr_\m
\ee
of the bundle $\op\w^r T^*Z\ot TZ\to Z$. Tangent-valued 1-forms are
usually termed the {\it (1,1) tensor fields}.

\begin{example} There is the 1:1 correspondence between the
tangent-valued 1-forms on $Z$ and the linear bundle morphisms
\beq
\f:TZ\to TZ,\qquad 
\f: T_zZ\ni v\mapsto v\rfloor\f(z)\in T_zZ \label{29} 
\eeq
over $Z$. In particular, the {\it canonical tangent-valued 1-form}
$\th_Z= dz^\la\ot \dr_\la$ defines 
the identity morphism of $TZ$. 
\end{example}

\subsection{Tangent and cotangent bundles of bundles}

The tangent bundle $TY\to Y$ of a bundle $\Y$ has the {\it vertical
tangent subbundle} $VY$ given by the coordinate condition $\dot
x^\la =0$. It is coordinatized by $(x^\la,y^i,\dot y^i)$ with respect to
the holonomic fibre bases $\{\dr_i\}$. Given a bundle morphism
$\Phi:Y\to Y'$, the restriction of the tangent morphism $T\Phi$ to $VY\subset
TY$ is the {\it vertical tangent morphism}
\be
V\Phi:VY\to VY', \qquad \dot y'^i\circ V\Phi =\dot y^j\dr_j\Phi^i.
\ee
We shall utilize the notation
\beq
\dr_V=\dot y^j\frac{\dr}{\dr y^j}. \label{m141}
\eeq

\begin{example}
If $\Y$ and an affine bundle modelled on the
vector bundle $\ol Y\to X$, there are the canonical isomorphisms
\be
 V\ol Y=\ol Y\op\xx_X \ol Y, \qquad  VY=Y\op\xx_X\ol Y. 
\ee
\end{example}

The {\it vertical cotangent bundle} $V^*Y\to Y$ of $\Y$ is defined to be the
vector bundle dual to the vertical tangent bundle $VY\to Y$. It is not
a subbundle of the cotangent bundle $T^*Y$.
We shall denote by $\{\ol dy^i\}$ the
fibre bases for $V^*Y$ which are dual to the fibre bases $\{\dr_i\}$ for
$VY$. 

We have the following exact sequences:
\bea
&& 0\to VY\op\har_Y TY\op\to_Y Y\op\times_X TX\to 0,
\label{1.8a} \\
&& 0\to Y\op\times_X T^*X\op\har_Y T^*Y\op\to_Y V^*Y\to 0.
\label{1.8b}
\eea
Any splitting
\bea
&& \G:TX\op\har_Y TY, \qquad
\dr_\la\mapsto\dr_\la +\G^i_\la (y)\dr_i, \label{m3a} \\
&& \G:V^*Y\op\har_Y T^*Y,\qquad
\ol dy^i\mapsto dy^i -\G^i_\la (y)dx^\la, \label{m3b}
\eea
of these sequences corresponds to the choice of a connection on the bundle
$Y\to X$.

Let us consider the bundles $TT^*X$ and $T^*TX$. Given the coordinates
$(x^\la, p_\la=\dot x_\la)$ of $T^*X$ 
and $(x^\la, v^\la=\dot x^\la)$ of $TX$,
these bundles are  coordinatized by $(x^\la, p_\la, \dot x^\la, \dot p_\la)$
and $(x^\la, v^\la,
\dot x_\la, \dot v_\la)$ respectively. By inspection of the
coordinate transformation laws, one can show that
they are isomorphic over $TX$ (see also \cite{cram,kijo}): 
\be
\bx{TT^*X\op=_{TX}T^*TX, \quad p_\la\llra\dot v_\la, \quad \dot p_\la
\llra\dot x_\la.} 
\ee

Given a bundle $\Y$, the similar isomorphism of the bundles $VV^*Y$ and
$V^*VY$ over $VY$ takes place. In 
coordinates $(x^\la, y^i, p_i=\dot y_i)$ of
$V^*Y$ and $(x^\la, y^i, v^i=\dot y^i)$ of $VY$, this isomorphism reads
\beq
\bx{VV^*Y\op=_{VY}V^*VY, \quad p_i\llra\dot v_i, 
\quad \dot p_i\llra\dot y_i.}
\label{m6}
\eeq

\subsection{Forms and vector fields on bundles}

Let $\pi: \Y$ be a bundle coordinatized by $(x^\la, y^i)$. 

A vector field $u$ on $\Y$ is termed {\it projectable} when there is a
vector field $\tau$ on $X$ such that 
\be
T\pi\circ u=\tau\circ \pi.
\ee
Its coordinate expression reads
\be
u=u^\m (x)\dr_\m +u^i(y)\dr_i, \qquad \tau=u^\mu(x)\dr_\mu.
\ee
A {\it vertical} vector field $u= u^i\dr_i$ is 
a projectable vector field over
the zero vector field on $X$.

We mention the following types of forms on a bundle $\Y$:
\begin{itemize}
\item {\it horizontal} forms
\be
\phi : Y\to\op\w^r T^*X, \quad
\phi =\phi_{\la_1\cdots\la_r}(y)dx^{\la_1}\w\cdots\w dx^{\la_r},
\ee
\item {\it tangent-valued horizontal} forms
\be
 &&\phi : Y\to\op\w^r T^*X\op\ot_Y TY,\\
 &&\phi =dx^{\la_1}\w\cdots\w dx^{\la_r}\ot
[\phi_{\la_1\cdots\la_r}^\m (y)\dr_\m +
\phi_{\la_1\cdots\la_r}^i(y) \dr_i],
\ee
\item {\it vertical-valued horizontal} forms
\be
\phi : Y\to\op\w^r T^*X\op\otimes_Y VY,\quad
\phi =\phi_{\la_1\cdots\la_r}^i(y)
dx^{\la_1}\w\cdots\w dx^{\la_r}\ot\dr_i,
\ee
\item {\it pullback-valued forms}
\ben
&&Y\to \op\w^r T^*X\op\ot_Y TX,  \qquad
 \phi =\phi_{\la_1\dots\la_r}^\m (y)
dx^{\la_1}\w\cdots\w dx^{\la_r}\otimes \dr_\m,\label{1.11}\\
&&Y\to \op\w^r T^*X\op\ot_Y V^*Y, \qquad
 \phi =\phi_{\la_1\dots\la_ri}(y)
dx^{\la_1}\w\cdots\w dx^{\la_r}\ot \ol dy^i.\label{87}
\een
\end{itemize}

Horizontal 1-forms on $Y$ are called {\it semi-basic forms}. If such a
form is the pullback $\pi^*\f$ of a 1-form $\f$ on $X$, it is said to be a
{\it basic form}. If there is no danger of confusion, we shall denote
the pullbacks
$\pi^*\f$ onto $Y$ of forms $\f$ on $X$ by the same symbol $\f$.

Vertical-valued horizontal 1-forms
\beq
\sigma : Y\to T^*X\op\otimes_Y VY, \quad \si = \si_\la^i dx^\la\otimes\dr_i,
\label{m23}
\eeq
 are termed the {\it soldering forms}. 

Note that the forms (\ref{1.11}) are not tangent-valued forms, and the forms 
(\ref{87}) are not exterior forms. The pullbacks
\be
\phi =\phi_{\la_1\cdots\la_r}^\mu (x)
dx^{\la_1}\wedge\cdots\wedge dx^{\la_r}\otimes\dr_\mu
\ee
of tangent-valued forms on $X$ onto $Y$ exemplify the pullback-valued
forms (\ref{1.11}). 

Horizontal forms of degree $n=\dim X$ on a 
bundle $Y\to X$ are called the {\it horizontal densities}.
In sequel, we shall exploit the notation
\beq
\om=dx^1\w \cdots\w dx^n, \qquad 
\om_\la =\dr_\la\rfloor\om, \qquad \dr_\m\rfloor\om_\la = \om_{\m\la}.
\label{m112}
\eeq

\subsection{Distributions}

Let $Z$ be an $m$-dimensional manifold.
A $k$-codimensional {\it smooth distribution}
$\cT$ on $Z$ is  defined to be a subbundle 
of rank $m-k$ of the tangent bundle
$TZ$. A smooth distribution
$\cT$ is called {\it involutive} if $[u,u']$ is a section of $\cT$
whenever $u$ and $u'$ are sections of $\cT$.

Let $\cT$ be a $k$-codimensional distribution on $Z$. 
Its annihilator $\cT^*$
is a
$k$-dimensional subbundle of $T^*Z$ 
called the {\it Pfaffian system}. It means
that, on a neighborhood $U$ of 
every point $z\in Z$, there exist $k$ linearly
independent sections $\f_1,\dots,\f_k$ of $\cT^*$ such that 
\be
\cT_z\mid_U=\op\cap_j\Ker\f_j. 
\ee
Let $\cC(\cT)$ be the ideal of $\cO(Z)$ generated by sections of $\cT^*$.

\begin{proposition}
A smooth distribution $\cT$ is involutive iff the ideal $\cC(\cT)$ is
{\it differential}, that is, $d\cC(\cT)\subset \cC(\cT)$ \cite{warn}.
\end{proposition}

\begin{remark} Given an involutive $k$-codimensional distribution $\cT$ on
$Z$, the quotient 
$TZ/\cT$ is a $k$-dimensional vector bundle called the {\it
transversal bundle} of $\cT$. There is the exact sequence
\beq
0\to \cT\har TZ\to TZ/\cT\to 0. \label{m105}
\eeq
Given a bundle $\Y$, its vertical tangent bundle $VY$ exemplifies an
involutive distribution on $Y$. In this case, 
the exact sequence (\ref{m105})
is just the exact sequence (\ref{1.8a}).
\end{remark}

 A submanifold $N$ of $Z$ is called the {\it integral
manifold} of a distribution
$\cT$ on $Z$ if the tangent spaces to $N$ coincide with the fibres of this
distribution at each point of $N$.

\begin{theorem}\label{to.1} Let $\cT$ be a smooth involutive distribution on
$Z$. For any point $z\in Z$, there exists a maximal integral manifold of 
$\cT$ passing through $z$ \cite{kamb,libe,warn}. 
\end{theorem}

In view of this fact, involutive distributions 
are also called {\it completely integrable distributions}.

\begin{corollary}\label{c11.0}  
Every point $z\in Z$ has an open neighborhood
$U$ which is a domain of a coordinate chart $(z^1,\dots,z^m)$ such that
the restrictions of $\cT$ and $\cT^*$ to 
$U$ are generated by the $m-k$ vector
fields $\frac{\dr}{\dr z^1}, \dots,\frac{\dr}{\dr z^{m-k}}$ and the $k$
Pfaffian forms $dz^{m-k+1},\dots, dz^m$ respectively.
\end{corollary}

In particular, it follows that integral manifolds of an involutive 
distribution constitute a foliation. 
Recall that a $k$-codimensional {\it foliation} 
on a $m$-dimensional manifold
$Z$ is  a partition of $Z$ into connected 
leaves $F_\iota$ with the following
property.  Every point of $Z$ has an open neighborhood $U$ which is a domain
of a coordinate chart
$(z^\al)$ such that, for every leaf $F_\iota$, the components
$F_\iota\cap U$ are described by the equations
$z^{m-k+1}=$const. ,..., $z^m=$const. \cite{kamb,rein}. 
Note that leaves of a foliation fail to be imbedded submanifolds in general.

\begin{example}
Every projection $\pi:Z\to X$ defines a
foliation whose leaves are the fibres $\pi^{-1}(x)$, $x\in X$. 
\end{example}

\begin{example}\label{line} Every nowhere vanishing 
vector field $u$ on a manifold $Z$ defines a 1-dimensional
involutive distribution on $Z$. Its integral manifolds are the integral
curves of $u$. In virtue of the
Corollary \ref{c11.0}, around each point $z\in Z$, there exist local
coordinates $(z^1,\dots,z^m)$ of a neighborhood of 
$z$ such that $u$ is given by $u=\frac{\dr}{\dr z^1}$.
\end{example}

\subsection{First order jet manifolds}

Differential operators, differential equations and Lagrangian formalism are
conventionally phrased in terms of jet manifolds
\cite{gols,kola,pomm,sard0,sard93,saun}.  

Given a bundle $\Y$, its first order {\it jet manifold} $J^1Y$ comprises the
equivalence classes $j^1_xs$, $x\in X$, of sections $s:X\to Y$ 
so that sections
$s$ and $s'$ belong to the same class iff
\be
Ts\mid _{T_xX} =Ts'\mid_{T_xX}.
\ee
Roughly speaking, sections $s,s'\in j^1_xs$  are identified by their values
$s^i(x)={s'}^i(x)$ and the values of their partial derivatives
$\dr_\mu s^i(x)=\dr_\mu{s'}^i(x)$
at the point $x$ of $X$. There are the natural fibrations
\be
\pi_1:J^1Y\ni j^1_xs\mapsto x\in X, \qquad 
\pi_{01}:J^1Y\ni j^1_xs\mapsto s(x)\in Y. 
\ee

Given bundle coordinates $(x^\la,y^i)$ of $Y$, the jet manifold $J^1Y$ 
is equipped with the adapted coordinates
\ben
&&(x^\la,y^i,y_\la^i), \quad (y^i,y_\la^i)(j^1_xs)=(s^i(x),\dr_\la
s^i(x)), \nonumber\\
&&{y'}^i_\la = \frac{\dr x^\m}{\dr{x'}^\la}(\dr_\m +y^j_\m\dr_j)y'^i.
\label{50}
\een
A glance at (\ref{50}) shows that $J^1Y\to Y$
is an affine bundle modelled on the vector bundle
\beq
T^*X \op\otimes_Y VY\to Y.\label{23}
\eeq

\begin{proposition}
There exist the canonical bundle monomorphisms 
\ben
&& \la:J^1Y\op\har_Y
T^*X \op\otimes_Y TY,\quad \la=dx^\la\ot d_\la=dx^\la
\otimes (\dr_\la + y^i_\la \dr_i),\label{18} \\
&& \th:J^1Y \op\har_Y T^*Y\op\otimes_Y VY,\quad
\th=\th^i \otimes \dr_i=(dy^i- y^i_\la dx^\la)\ot\dr_i,\label{24}
\een
where $d_\la$ are called the {\it total derivatives}.
\end{proposition}

\begin{remark} The total derivatives obey the relations
\be
d_\la\circ d= d\circ d_\la, \qquad d_\la(\f\w\si) =d_\la\f\w\si +\f\w
d_\la\si, \qquad \f,\si\in \cO(Y).
\ee
\end{remark}

The monomorphisms (\ref{18}) and (\ref{24}) enable us to express jets into
the tangent-valued forms.

Let $\Phi: Y\to Y'$ be a bundle morphism over a diffeomorphism $f$ of $X$.
The {\it jet prolongation} of $\Phi$ is the morphism
\be
J^1\Phi : J^1Y  \to J^1Y',\quad
J^1\Phi :j_x^1s\mapsto j_{f(x)}^1(\Phi\circ s\circ f^{-1}).
\ee
Its coordinate expression is given by (\ref{50}).

Every section $s$ of a bundle $Y\to X$ admits the
jet prolongation to the section  
\be
(J^1s)(x) = j_x^1s, \qquad (y^i,y_\la^i)\circ J^1s= (s^i(x),\dr_\la s^i(x)),
\ee
of the bundle $J^1Y\to X$. A section of $J^1Y\to X$ is called
{\it holonomic} if it is the jet prolongation of a section of $\Y$.

Every projectable vector field $u = u^\la\dr_\la + u^i\dr_i$
on $\Y$ has the {\it jet lift} to the vector field
\ben
&&\ol u =r\circ J^1 u: J^1Y\to J^1TY\to TJ^1Y,\nonumber \\
&& \bx{\ol u =
u^\la\dr_\la + u^i\dr_i + (d_\la u^i
 - y_\m^i\dr_\la u^\m)\dr_i^\la,} \label{1.21}
\een
on the bundle $J^1Y\to X$.
In the definition of $\ol u$, we use the bundle morphism
\be
r: J^1TY\to TJ^1Y, \qquad \dot y^i_\la\circ r_1 = (\dot y^i)_\la-y^i_\m\dot
x^\m_\la.
\ee
In particular, there exists the canonical isomorphism
\beq
VJ^1Y=J^1VY, \qquad \dot y^i_\la=(\dot y^i)_\la.\label{1.22}
\eeq
As a consequence, the jet lift (\ref{1.21}) of a vertical vector field $u$
on $Y\to X$ coincides with its jet prolongation 
\be
\ol u=J^1 u=u^i\dr_i +d_\la u^i\dr^\la_i.
\ee

If a bundle $Y\to X$ is endowed with an algebraic structure, the jet bundle
$J^1Y\to X$ inherits this algebraic structure  due to the jet
prolongations of the corresponding morphisms. For instance, 
if $Y$ is a vector bundle, $J^1Y\to X$ does as well. If $Y$ is an affine
bundle modelled on a vector bundle $\ol Y\to X$, then
$J^1Y\to X$ is an affine bundle modelled on 
the vector bundle $J^1\ol Y\to X$.

The canonical monomorphisms (\ref{18}) and (\ref{24}) determine the {\it
canonical splitting} 
\beq
\bx{\dot x^\la\dr_\la
+\dot y^i\dr_i =\dot x^\la(\dr_\la +y^i_\la\dr_i) + (\dot y^i-\dot x^\la
y^i_\la)\dr_i} \label{1.20}
\eeq
of the pullback $J^1Y\op\times_Y TY$ and the dual splitting 
\beq
\bx{\dot x_\la dx^\la
+\dot y_i dy^i =(\dot x_\la + \dot y_iy^i_\la)dx^\la + \dot y_i(dy^i-
y^i_\la dx^\la)} \label{34}
\eeq
of the pullback $J^1Y\op\times_Y T^*Y$.
In particular, we get the canonical splitting of a vector field on $Y$:
\beq
u =u^\la\dr_\la +u^i\dr_i=u_H +u_V =u^\la (\dr_\la +y^i_\la
\dr_i)+(u^i - u^\la y^i_\la)\dr_i. \label{31}
\eeq

\subsection{Second order jet manifolds}

The {\it repeated jet manifold}
$J^1J^1Y$  is defined to be the jet manifold of the bundle
$J^1Y\to X$. It 
is provided with the adapted coordinates $(x^\la ,y^i,y^i_\la
,y_{(\m)}^i,y^i_{\la\m})$.

There are the projections
\ben
&&\pi_{11}:J^1J^1Y\to J^1Y, \qquad y_\la^i\circ\p_{11} = y_\la^i,
\label{S1}\\
&&J^1\pi_{01}:J^1J^1Y\to J^1Y,\qquad
y_\la^i\circ J^1\pi_{01} = y_{(\la)}^i. \label{S'1}
\een
They coincide on the
{\it sesquiholonomic subbundle} $\wh J^2Y\to J^1Y$ of $J^1J^1Y$ which is 
given by the coordinate conditions $y^i_{(\la)}= y^i_\la$. 
It is coordinatized by $(x^\la ,y^i, y^i_\la,y^i_{\la\m})$.

The {\it second order jet manifold} $J^2Y$ 
of a bundle $\Y$ is the subbundle of $\wh J^2Y\to J^1Y$ defined
by the coordinate conditions 
$y^i_{\la\m}=y^i_{\m\la}$. It is coordinatized by
$(x^\la ,y^i, y^i_\la,y^i_{\la\leq\m})$
together with the transition functions
\be
{y'}_{\la\m}^i= \frac{\dr x^\al}{\dr{x'}^\m}(\dr_\al +y^j_\al\dr_j
+y^j_{\nu\al}\dr^\nu_j){y'}^i_\la. 
\ee
The second order jet manifold $J^2Y$  of $Y$ comprises
the equivalence classes  $j_x^2s$ of sections $s$ of $Y\to X$ such that
\be
y^i_\la (j_x^2s)=\dr_\la s^i(x),\qquad
y^i_{\la\m}(j_x^2s)=\dr_\m\dr_\la s^i(x).
\ee
In other words, two sections $s,s'\in j^2_xs$ are identified by their values
and the values of their first and second order 
derivatives at the point $x\in X$.

Let $\Phi:Y\to Y'$ be a bundle morphism over
a diffeomorphism of $X$ and $J^1\Phi$ its jet prolongation.
Let us consider the jet prolongation $J^1J^1\Phi:J^1J^1Y \to
J^1J^1Y'$ of $J^1\Phi$.
Restricted to the second order jet manifold
$J^2Y$, the morphism $J^1J^1\Phi$ takes its values in $J^2Y'$.
It is called the {\it second order jet prolongation} $J^2\Phi$ of $\Phi$.

Similarly, the repeated jet prolongation
$J^1J^1s$ of a section $s$ of $Y\to X$ is a  section of the bundle
$J^1J^1Y\to X$. It takes its values into $J^2Y$ and
defines the following second order jet prolongation of $s$:
\be
(J^2s)(x)=(J^1J^1s)(x)=j^2_xs.
\ee

\subsection{Ehresmann connections}

A {\it connection} $\G$ on $Y$ is usually defined to be a splitting
(\ref{m3a}) (or (\ref{m3b})) of the exact sequence  (\ref{1.8a}) (or
(\ref{1.8b})). There are the corresponding decompositions 
\beq
TY=\G(TX)\op\oplus_Y VY, \qquad T^*Y=T^*X\op\oplus_Y \G(V^*X).\label{11}
\eeq

It is readily observed that
the canonical splittings (\ref{1.20}) -- (\ref{34}) of 
$TY$ and $T^*Y$ over the jet bundle
$J^1Y\to Y$ enable us to recover the splittings (\ref{11}) 
by means of a section 
\beq
\bx{\G =dx^\la\otimes(\dr_\la +\G^i_\la (y)\dr_i),} \qquad 
{\G'}^i_\la = (\frac{\dr{y'}^i}{\dr y^j}\G_\m^j +
\frac{\dr{y'}^i}{\dr x^\m})\frac{\dr x^\m}{\dr{x'}^\la}, \label{37}
\eeq
of this jet bundle.
Substituting $\G$ (\ref{37}) into  (\ref{1.20}) -- (\ref{34}), we obtain the
familiar splittings
\ben
&&\dot x^\la\dr_\la +\dot y^i\dr_i = \dot x^\la (\dr_\la
+\G^i_\la\dr_i) + (\dot y^i-\dot x^\la\G^i_\la)\dr_i, \nonumber\\
&&\dot x_\la dx^\la +\dot y_idy^i = (\dot x_\la +\G^i_\la\dot
y_i)dx^\la + \dot y_i(dy^i-\G^i_\la dx^\la) \label{9}
\een
corresponding to (\ref{11}).

Hereafter, we follow the notion of a connection on $\Y$ as a section of
the jet bundle $J^1Y\to Y$. It is called the {\it Ehresmann connection}.

\begin{example}\label{c1}
Let $Y\to X$ be a vector bundle. A linear connection on $Y$ reads
\beq
\G=dx^\la\otimes[\dr_\la-\G^i{}_{j\la}(x)y^j\dr_i]. \label{8}
\eeq
Let $Y\to X$ be an affine bundle modelled on a vector bundle
$\ol Y\to X$. An affine connection on $Y$ reads
\be
\G=dx^\la\otimes[\dr_\la+(-\G^i{}_{j\la}(x)y^j+\G^i{}_\la (x)) \dr_i],
\ee
where $\ol\G=dx^\la\otimes[\dr_\la-\G^i{}_{j\la}(x)\ol y^j\dr_i]$
is a linear connection on $\ol Y$.
\end{example}

Since the affine jet bundle $J^1Y\to Y$ is modelled on the vector bundle
(\ref{23}), Ehresmann connections on $Y\to X$ constitute an affine space
modelled on the linear space of soldering forms on $Y$. If
$\G$ is a connection and $\si$ 
is a soldering form (\ref{m23}) on $Y$, its sum
\be
\G+\si=dx^\la\otimes[\dr_\la+(\G^i_\la +\si^i_\la)\dr_i]
\ee
is a connection on $Y$. Conversely, if $\G$ and $\G'$ are
connections on $Y$, then
\be
\G-\G'=(\G^i_\la -{\G'}^i_\la)dx^\la\otimes\dr_i
\ee
is a soldering form.

We mention the following operations with connections. 

(i) Let $\G$ be a connection
on $Y\to X$ and $\G'$ a connection on $Y'\to X$. There
exists the {\it product connection} $\G\times\G'$ on $Y\op\times_X Y$.

(ii) Every linear connection $\G$ on a vector bundle $Y\to X$
yields the {\it dual linear connection} 
\be
\G^*_{i\la}=\G^j{}_{i\la}(x)y_j
\ee
on the dual vector bundle $Y^*\to X$.

\begin{example}\label{c4} A linear 
connection $K$ on the tangent bundle $TX$
of a manifold $X$ and the dual connection $K^*$ to $K$ on the cotangent 
bundle $T^*X$ read
\beq
\bx{K^\al_\la=-K^\al{}_{\nu\la}(x)\dot x^\nu,\qquad
K^*_{\al\la}=K^\nu{}_{\al\la}(x)\dot x_\nu.} \label{408}
\eeq
\end{example}

(iii) Given a connection $\G$ on $\Y$, the vertical tangent
morphism $V\G$ yields the {\it vertical connection}
\beq
 \bx{V\G =dx^\la\otimes(\dr_\la
+\G^i_\la\frac{\dr}{\dr y^i}+\dr_V\G^i_\la 
\frac{\dr}{\dr \dot y^i}),} \qquad
\dr_V\G^i_\la= \dot y^j\dr_j\G^i_\la,\label{43}
 \eeq
on the bundle $VY\to X$ due to the canonical isomorphism
(\ref{1.22}). The dual {\it covertical connection} on the bundle 
$V^*Y\to X$ reads
\beq
\bx{V^*\G =dx^\la\otimes(\dr_\la +\G^i_\la\frac{\dr}{\dr
y^i}-\dot y_i\dr_j\G^i_\la \frac{\dr}{\dr \dot y_j}).} \label{44}
\eeq

(iv) For every connection $\G$ on $Y\to X$,
one can construct its {\it jet lift} $J\G$ onto the bundle $J^1Y\to X$ as
follows. Note that 
the jet prolongation $J^1\G$ of the connection $\G$ on $Y$ is
 a section of the repeated jet bundle (\ref{S'1}), but not of the
bundle $\pi_{11}$ (\ref{S1}).
Let $K^*$ be a linear symmetric connection (\ref{408}) on the cotangent 
bundle $T^*X$ of $X$.
There exists the  affine bundle morphism
\be
&&r_K: J^1J^1Y\to J^1J^1Y, \qquad  r_K\circ r_K=\Id_{J^1J^1Y},\\
&&(y^i_\la ,y_{(\m)}^i,y^i_{\la\m})\circ r_K=
(y^i_{(\la)} ,y_\m^i,y^i_{\m\la}+ K^\al{}_{\la\m}(y^i_\al - y^i_{(\al)})).
\ee
We set
\beq
\bx{J\G=r_K\circ J^1\G= dx^\m\ot
[\dr_\mu+\G^i_\mu\dr_i +(d_\la\G^i_\m -
K^\al{}_{\la\mu} (y^i_\al-\G^i_\al)) \dr_i^\la].} \label{59}
\eeq

The {\it curvature} of a connection $\G$ is given by the
horizontal vertical-valued 2-form
\ben
&&R  =\frac12\sum R^i_{\la\m} dx^\la\w dx^\m\otimes\dr_i,\nonumber \\
&& \bx{R^i_{\la\m}= \dr_\la\G^i_\m -\dr_\m\G^i_\la +\G^j_\la\dr_j\G^i_\m
-\G^j_\m\dr_j\G^i_\la.} \label{13}
\een
In particular, the curvature of the linear connection (\ref{8}) reads
\be
&& R^i_{\la\m}(y)=-R^i{}_{j\la\m}(x)y^j, \\
&&R^i{}_{j\la\m}=\dr_\la\G^i{}_{j\m} -\dr_\m\G^i{}_{j\la}
+\G^k{}_{j\mu}\G^i{}_{k\la} -\G^k{}_{j\la}\G^i{}_{k\mu}.
\ee

A connection $\G$ on $Y\to X$ yields the first order differential operator
\beq
D_\G:J^1Y\to T^*X\op\otimes_Y VY, \qquad \bx{D_\G =(y^i_\la
-\G^i_\la)dx^\la\otimes\dr_i,}\label{38}
\eeq
called the {\it covariant differential} relative to $\G$.
The corresponding {\it covariant
derivative} of a section $s$ of  $Y$ is
\be
\nabla_\G s=D_\G\circ J^1s=[\dr_\la s^i-
(\G\circ s)^i_\la]dx^\la\otimes\dr_i. 
\ee
A local section $s$ of a $Y\to X$
is said to be an {\it integral section} for a
connection $\G$ on $Y$ if  $\G\circ s=J^1s$, that is, $\nabla_\G s=0$.

\begin{remark}\label{coneq}
Every connection $\G$ on the bundle $Y\to X$ defines a system of first order
differential equations on $Y$ (in the spirit of \cite{gols,kras,pomm}) 
which is an imbedded subbundle
$\G(Y)=\Ker D_\G$ of the jet bundle $J^1Y\to Y$. It is given by the 
coordinate relations
\beq
y^i_\la =\G^i(y). \label{39}
\eeq
Integral sections for $\G$ are local solutions of (\ref{39}), and
{\it vice versa}.
\end{remark}

\subsection{Curvature-free connections}

Every connection $\G$ on $\Y$, by definition, yields the {\it horizontal
distribution} $\G(TX)\subset TY$ (\ref{m3a}) on $Y$. It is generated
by horizontal lifts
\be
\tau_\G= \tau^\la(\dr_\la +\G^i_\la\dr_i) 
\ee
onto $Y$ of vector fields $\tau=\tau^\la\dr_\la$ on
$X$. The associated Pfaffian system is locally generated by the forms
$(dy^i-\G^i_\la dx^\la)$.

\begin{proposition}\label{flat}
The horizontal distribution $\G(TX)$
is involutive iff $\G$ is a curvature-free connection.
\end{proposition}

\begin{proof}
Straighforward calculations show that
\be
[\tau_\G, {\tau'}_\G] = ([\tau,\tau'])_\G
\ee
iff the curvature $R$ (\ref{13}) of $\G$ vanishes everywhere.
\end{proof}

\begin{remark}
Obviously, not every bundle admits a curvature-free connection.
If a principal bundle over a simply connected base (i.e., its
first homotopy group is trivial) admits a curvature-free connection, this
bundle is trivializable \cite{koba}.
\end{remark}

In virtue of Theorem \ref{to.1}, the horizontal distribution defined by a
curvature-free connection is completely integrable.
The corresponding foliation on $Y$ is transversal to the
foliation defined by the fibration $\pi:\Y$. 
It is called the {\it horizontal
foliation}. Its leaf through a point
$y\in Y$ is defined locally by the integral section $s_y$ of the connection
$\G$ through $y$. 
Conversely, let $Y$ admits a horizontal foliation such that, for each
point $y\in Y$, the leaf of this foliation through $y$ is locally defined 
by some section $s_y$ of $Y\to X$ through $y$. Then, the map
\be
\G:Y\to J^1Y, \qquad \G(y)=j^1_ss_y, \qquad \pi(y)=x. 
\ee
is well defined. This is a curvature-free connection on $Y$. 

\begin{corollary}\label{horfol} There is the 1:1 correspondence between the
curvature-free connections and the horizontal foliations on a bundle $\Y$.
\end{corollary}

Given a horizontal foliation on $\Y$, there exists the  associated atlas of
bundle coordinates $(x^\la, y^i)$  of $Y$ such that (i) every leaf of this
foliation is locall generated by the equations $y^i=$const. and (ii) the
transition functions
$y^i\to {y'}^i(y^j)$ are independent on the coordinates $x^\la$ of the base
$X$ \cite{kamb}. It is called the {\it atlas of constant local
trivializations}.  Two such atlases are said to be equivalent if their union
also is an atlas of  constant local trivializations. They are
associated with the same horizontal foliation.

\begin{corollary}\label{can}
There is the 1:1 correspondence between 
the curvature-free connections $\G$ on
a bundle $Y\to X$ and the equivalence classes 
of atlases $\Psi_c$ of constant
local trivializations of
$Y$ such that $\G^i_\la=0$ relative to the
coordinates of the corresponding atlas $\Psi_c$ \cite{cana}.
\end{corollary}

Connections on a bundle over a 1-dimensional base $X^1$ are curvature-free
connections.

\begin{example}\label{rcon}
Let $Y\to X^1$ be such a bundle ($X^1=\R$ or $X^1=S^1$, see Remark
\ref{onedim}). It is coordinatized by $(t,y^i)$, where $t$ is either
the canonical parameter of $\R$ or the standard local coordinate of $S^1$
together with the transition functions
$t'=t+$const. Relative to this coordinate, the base $X^1$ is
provided with the standard 
vector field $\dr_t$ and the standard 1-form $dt$. 
Let $\G$  be a connection
on $Y\to X^1$. In virtue of Proposition \ref{flat}, such a connection
defines a horizontal foliation on
$Y\to X^1$. Its leaves are the integral curves of the horizontal lift
\beq
\tau_\G=\dr_t + \G^i\dr_i \label{m134}
\eeq
of $\dr_t$ by $\G$ (see Example \ref{line}). The corresponding Pfaffian
system is locally generated by the forms  $(dy^i-\G^i dt)$. There exists an
atlas of constant local trivializations $(t,y^i)$ such that $\G^i=0$ and
$\tau_\G=\dr_t$ relative to these coordinates. 
\end{example}

A connection
$\G$ on $Y\to X^1$ is called {\it complete} if the horizontal vector field
(\ref{m134}) is complete. 

\begin{proposition}\label{complcon}
Every trivialization of  $Y\to \R$ defines a complete connection.
Conversely, every complete connection on $Y\to
\R$ defines a trivialization $Y\simeq \R\xx M$. 
The vector field (\ref{m134})
comes to the vector field $\dr_t$ on $\R\xx M$.
\end{proposition}

\begin{proof}
Every trivialization of $Y\to\R$ defines a one-parameter group of
isomorphisms of $Y\to\R$ over $\Id_\R$, and hence a complete connection.
Conversely, let $\G$ be a complete connection on $Y\to\R$. The vector field
$\tau_\G$ (\ref{m134}) is the generator of a 1-parameter group $G_\G$ which
acts freely on $Y$. The orbits of this action are of course the integral
sections of $\tau_\G$. Hence we get a projection $Y\to M=Y/G_\G$ which,
together with the projection $Y\to\R$, defines a trivialization 
$Y\simeq \R\xx M$.
\end{proof}

\subsection{Composite connections}

Let us consider a bundle $\pi:\Y$ which admits a {\it composite fibration}
\beq
 Y\lra^{\pi_{Y\Si}} \Si\lra^{\pi_{\Si X}} X, \label{1.34}
\eeq
where $Y\to\Si$ and $\Si\to X$ are bundles. It is equipped with the bundle
coordinates $(x^\la,\si^m,y^i)$ together with the transition functions
\be
x^\la\to {x'}^\la(x^\mu), \qquad
\si^m\to {\si'}^m(x^\mu,\si^n), \qquad
y^i\to {y'}^i(x^\mu,\si^n,y^j),
\ee
 where $(x^\m,\si^m)$ are bundle coordinates  of $\Si\to X$. 

\begin{example}
We have the composite bundles
\be
TY\to Y\to X, \qquad VY\to Y\to X, \qquad J^1Y\to Y\to X.
\ee
\end{example}

Let
\beq
A=dx^\la\ot(\dr_\la+ A^i_\la\dr_i)
+ d\si^m\ot(\dr_m+A^i_m\dr_i) \label{Q1}
\eeq
be a connection on the bundle $Y\to \Si$ and
\be
 \G=dx^\la\ot(\dr_\la + \G^m_\la\dr_m)
\ee
a connection on the bundle $\Si\to X$.
Given a vector field $\tau$ on $X$, let us consider its horizontal lift
$\tau_\G$ onto $\Si$ by $\G$ and then the horizontal lift
$(\tau_\G)_A$ of $\tau_\G$ onto $Y$ by the connection (\ref{Q1}).

\begin{proposition} There exists the connection 
\beq
\bx{\g=dx^\la\ot[\dr_\la+\G^m_\la
\dr_m + (A^i_m\G^m_\la + A^i_\la)
\dr_i].} \label{1.39}
\eeq
 on $\Y$ such that the horizontal lift $\tau_\g$ onto $Y$ of
any vector field $\tau$ on $X$ 
consists with the above mentioned lift $(\tau_\G)_A$
\cite{sard93,sard95}. It is called the {\it composite connection}.
\end{proposition}

Given a composite bundle $Y$ (\ref{1.34}), the exact sequence
\be
0\to VY_\Si\har VY\to Y\op\times_\Si V\Si\to 0
\ee
over $Y$ take place, where $VY_\Si$ is the vertical tangent bundle of
$Y\to\Si$. Every connection (\ref{Q1}) on the bundle $Y\to\Si$ yields
the splitting
\be
&&VY=VY_\Si\op\oplus_Y (Y\op\times_\Si V\Si),\\
&&\dot y^i\dr_i + \dot\si^m\dr_m=
(\dot y^i -A^i_m\dot\si^m)\dr_i + \dot\si^m(\dr_m+A^i_m\dr_i). 
\ee
Due to this splitting, one can construct
the first order differential operator 
\ben
&&\wt D={\rm pr}_1\circ D_\g
: J^1Y\to T^*X\op\otimes_Y VY \to T^*X\op\otimes_Y VY_\Si,\nonumber\\
 &&\bx{\wt D=   dx^\la\ot (y^i_\la-
A^i_\la -A^i_m\si^m_\la)\dr_i,} \label{7.10}
 \een
on the composite manifold $Y$,
where $D_\g$ is the covariant differential (\ref{38}) relative to the
composite connection (\ref{1.39}). We call $\wt D$ the {\it vertical
covariant differential}.

\newpage

\section{Symplectic geometry}

This Section aims to recall the basic notions of 
symplectic geometry which we
shall need in sequel \cite{abra,arno,libe,vais}.

\subsection{Jacobi structure}

Let $Z$ be a manifold. The {\it Jacobi bracket} (or the {\it Jacobi
structure}) on $Z$ is defined to be a bilinear map 
\be
S(Z)\xx S(Z)\ni (f,g)\to \{f,g\}\in S(Z),
\ee
where $S(Z)$ is the linear space of smooth
functions on $Z$, which satisfies the following conditions:

(A1) $\{g,f\}=-\{f,g\}$ (skew-symmetry),

(A2) $\{f,\{g,h\}\} + \{g,\{h,f\}\} +\{h,\{f,g\}\}=0$ (Jacobi identity),

(A3) the support of $\{f,g\}$ is contained in the intersection of the
supports of $f$ and $g$.

\begin{proposition}\label{p11.1} Every Jacobi bracket on a manifold $Z$ is
uniquely defined in accordance with the relation 
\beq
\bx{\{f,g\}=w(df,dg) + u\rfloor(fdg-gdf)} \label{m50}
\eeq
by a vector field $u$ and a bivector field
$w$ on $Z$ such that
\beq
\bL_uw=0, \qquad [w,w]=2u\w w \label{m76}
\eeq
\cite{kiri76,lich78,marl}.
\end{proposition}

Taking $w=0$, every vector field $u$ on a 
manifold
$Z$ defines the Jacobi bracket (\ref{m50}). The relations (\ref{m76})
are obviously satisfied.

The Jacobi bracket (\ref{m50}) with $u=0$ is said to be a {\it Poisson
bracket} (see Section 3.3). Contact forms on an odd-dimensional manifold
generate Jacobi brackets which are not the Poisson ones (see next Section).

\subsection{Contact forms}

\begin{definition}
Given a $(2m+1)$-dimensional manifold $Z$, a {\it contact form} on
$Z$ is defined to be a Pfaffian form $\th$ such that the form
$\th\w(d\th)^m\neq$ everywhere on $Z$. The pair $(Z,\th)$ is called the
{\it contact manifold}.
\end{definition}

Note that a manifold $Z$ equipped with 
a contact form $\th$ is orientable, and
$\th\w(d\th)^m$ is a volume element. 

The assertion below is a variant of the well-known Darboux's theorem
\cite{libe}.

\begin{theorem}\label{darb} 
Let $(Z,\th)$ be a  $(2m+1)$-dimensional 
contact manifold. Every point $z$ of
$Z$ has an open neighborhood $U$ which is the domain of a coordinate chart 
$(z^0,\dots, z^{2m})$ such that the contact form $\th$ has the local
expression
\be
\th =dz^0 - \op\sum^m_{i=1} z^{m+i} dz^i 
\ee
on $U$. These coordinates are called {\it Darboux's coordinates}.
\end{theorem} 

If $\th$ is a contact form, its differential $d\th$ is a presymplectic form
of rank $2m$ (see Definition \ref{d11.1}).

\begin{proposition}
Let $\th$ be a contact form on $Z$. 
There exists the unique nowhere vanishing
vector field $E$ on $Z$ such that
\be
\bx{E\rfloor\th =1, \qquad E\rfloor d\th=0.} 
\ee
It is called the {\it Reeb vector field} of $\th$ \cite{libe}.
\end{proposition}

Relative to Darboux's coordinates, the Reeb vector field reads $E=\dr_0$.

\begin{proposition}\label{p11.2}
Every contact form $\th$ on an odd-dimensional manifold $Z$
yields the associated Jacobi structure on $Z$. It is defined by the
Reeb vector field $E$ of $\th$ and by the bivector field $w$ such that
\beq
w^\sh\f\rfloor\th=0, \qquad w^\sh\f\rfloor d\th =-(\f- (E\rfloor\f)\th)
\label{m55}
\eeq
for every $\f\in\cO^1(Z)$  \cite{marl}.
\end{proposition}

Relative to Darboux's coordinates, the Jacobi structure (\ref{m55}) reads
\be
\{f,g\}=\op\sum^m_{i=1}(\dr_{m+i}g\dr_if - \dr_{m+i}f\dr_ig)  + (\wt
g\dr_0f - \wt f\dr_0g), 
\ee
where
\be
\wt f=\op\sum^m_{i=1} z^{m+i}\dr_{m+i}f+f, 
\qquad \wt g=\op\sum^m_{i=1} z^{m+i}\dr_{m+i}g+g.
\ee

\subsection{Poisson structure}

According to (\ref{m76}), a bivector field $w$ on a manifold
$Z$ provides a {\it Poisson structure} if it obeys the condition
\be
[w,w]=0, \qquad w^{\m\la_1}\dr_\m w^{\la_2\la_3} +w^{\m\la_2}\dr_\m
w^{\la_3\la_1} +w^{\m\la_3}\dr_\m w^{\la_1\la_2} =0. 
\ee
It is called the {\it Poisson bivector} 
(see Example \ref{biv}). A manifold $Z$
equipped with a Poisson bivector $w$ is called the {\it Poisson manifold}
$(Z,w)$. 

Besides the conditions (A1 -- A3), the {\it Poisson bracket}
\be
\bx{\{f,g\}=w(df,dg)} 
\ee
satisfies also the Leibniz rule
\beq
\{h,fg\}=\{h,f\}g +f\{h,g\}. \label{m82}
\eeq

A Poisson structure defined by a Poisson bivector $w$ is said to be
{\it regular} if  the associated
 morphism $w^\sh:T^*Z\to TZ$ (\ref{m51}) has a
constant rank. Hereafter, by a Poisson structure is meant a regular Poisson
structure.

A Poisson structure is called {\it nondegenerate} if $w^\sh$ is an
isomorphism. If the Poisson structure $w$ is
nondegenerate, it induces the symplectic form $\Om$ on $Z$ defined by the
coordinate relation $\Om_{\al\bt}= (w^{-1})_{\bt\al}$, and {\it vice versa}
(see Proposition \ref{p11.5}). A nondegenerate Poisson
structure can exist only on an even-dimensional manifold (see next Section).

Note that there are no pullback or push-forward operations of Poisson
structures by manifold morphisms in general. The following assertion
deals with Poisson projections, whereas Theorem \ref{t11.1} is concerned
with Poisson injections.

\begin{proposition}\label{p11.3} Let $(Z,w)$ be a Poisson manifold and
$\pi:Z\to Y$ a projection. The following properties are equivalent:

(i) for every pair $(f,g)$ of functions 
on $Y$ and for each point $y\in Y$, the
restriction of the function $\{f\circ\pi,g\circ\pi\}$ to the fibre
$\pi^{-1}(y)$ is constant;

(ii) there exists a Poisson structure on $Y$ for which $\pi$ is a Poisson
morphism.

\noindent
If such a Poisson structure exists, it is unique \cite{libe}.
\end{proposition}

\begin{definition}\label{d11.3}
Given a function $f$ on a Poisson manifold $(Z,w)$, the image 
\be
\bx{\vt_f=w^\sh df, \qquad \vt_f= w^{\m\nu}\dr_\m f\dr_\nu,}
\ee
of its differential $df$ by the morphism $w^\sh$ is called the {\it
Hamiltonian vector field} of $f$.
\end{definition}

The Hamiltonian vector field $\vt_f$, by definition, obeys the relation
\beq
\bx{\vt_f\rfloor dg=\{f,g\}} \label{m80}
\eeq
for any function $g$ on $Z$. It is easy to see that
\beq
[\vt_f,\vt_g]=\vt_{\{f,g\}}. \label{m81}
\eeq
This relation provides the set of Hamiltonian 
vector fields with a Lie algebra
structure. Using (\ref{m82}) and (\ref{m81}), one can show that
\be
(\bL_{\vt_h}w)(df,dg) =\vt_h\rfloor d\{f,g\} -\{\vt_h\rfloor df, g\}
-\{f,\vt_h\rfloor dg\}=0.
\ee
It follows that a Hamiltonian vector field  generates a (local) 1-parameter
group of isomorphisms of the Poisson manifold $(Z,w)$.

The values of all Hamiltonian vector fields at all points of $Z$ constitute
the  {\it characteristic distribution of the Poisson manifold} $(Z,w)$. In
virtue  of the relation (\ref{m81}), 
this distribution is involutive. We have 
the following theorem.

\begin{theorem}\label{t11.1} 
The characteristic distribution of a Poisson manifold $(Z,w)$ is completely
integrable. The Poisson structure induces 
the symplectic structures on leaves
of the corresponding foliation of $Z$ \cite{vais}. It is called the
{\it symplectic foliation}.
\end{theorem}

Of course, the symplectic foliation has the adapted coordinates described
in Corollary \ref{c11.0}. Moreover, one can choose these coordinates in such
a way to bring the Poisson bracket in the following canonical form
\cite{vais,wein}.

\begin{proposition}\label{canpoiss} For any point $z$ of a Poisson manifold,
there are local coordinates $(z^1,\dots,z^r,y^1,\dots,y^k,p_1,\dots,p_k)$
near $z$ such that
\beq
\{y^i,y^j\}=\{p_i,p_j\}=\{y^i,z^\al\}=\{p_i,z^\al\}=\{z^\bt,z^\al\}=0, \quad 
\{p_i,y^j\}=\dl^j_i. \label{m111}
\eeq
\end{proposition}

\subsection{Symplectic structure}

\begin{definition}\label{d11.1}
A 2-form $\Om$ on a manifold $Z$ is called  {\it presymplectic}
if it is closed. A presymplectic form $\Om$ is said to be {\it symplectic} if
it is nondegenerate (see Example \ref{sympl}). 
\end{definition}  

A manifold $Z$ equipped with a symplectic 
[presymplectic] form is said to be a
{\it symplectic [ presymplectic] manifold}.

\begin{proposition}\label{p11.5} On every 
even-dimensional manifold $Z$, there
is the 1:1 correspondence between the symplectic forms $\Om$ and the Poisson
bivectors $w$ in accordance with the equalities
\be
\bx{w(\f,\si)=\Om(w^\sh\f,w^\sh\si), \quad
\Om(\vt,\nu)=w(\Om^\fl\vt,\Om^\fl\nu),} \quad \f,\si\in\cO^1(Z),
\quad \vt,\nu\in \cV^1(Z), 
\ee
(see relations (\ref{m52}) and (\ref{m51})) \cite{lich77}. 
\end{proposition}

In particular, the notion of a Hamiltonian vector field  may also be
introduced on symplectic manifolds.

\begin{definition}\label{d11.2} A vector field $\vt$ 
on a symplectic manifold
$(Z,\Om)$ is said to be {\it locally Hamiltonian} [{\it Hamiltonian}] if
the form $\vt\rfloor\Om $ is closed [exact].
\end{definition}

As an immediate consequence of this definition, we observe that:

(i) a vector field
$\vt$ is locally Hamiltonian iff it is an infinitesimal symplectomorphism,
that is, 
\be
\bx{\bL_\vt\Om= d(\vt\rfloor\Om)=0;}
\ee

(ii) a vector field $\vt$ is Hamiltonian iff 
it is a Hamiltonian vector field
in accordance with Definition \ref{d11.3}, i.e. $\vt=\vt_f$, where 
\be
\bx{df=-\vt_f\rfloor\Om, \qquad \vt_f=\dr^if\dr_i -\dr_if\dr^i.}
\ee

Note that Definition \ref{d11.2} of locally Hamiltonian and Hamiltonian
vector fields applies also to presymplectic manifolds.

\begin{example}\label{e11.0} Let $M$ be a manifold coordinatized by $(y^i)$
and let $T^*M$ be its cotangent 
bundle provided with the holonomic coordinates
$(y^i, p_i=\dot y_i)$. The cotangent bundle $T^*M$ is a well-known 
symplectic manifold equipped with the {\it canonical symplectic form}
\beq
\bx{\Om= dp_i\w dy^i} \label{m83}
\eeq
which is the differential of the canonical Liouville form
\beq
\bx{\th= p_i dy^i} \label{m85}
\eeq
on $T^*M$. 
These forms are canonical in the sense that the expressions
(\ref{m83}) and (\ref{m85}) are maintained 
under arbitrary transformations of
the coordinates
$y^i$ and the corresponding holonomic transformations of the coordinates
$p_i$. Furthermore, for every closed 2-form $\f$ on $M$, the form
$\Om +\f$ is also a symplectic form on $T^*M$. 
\end{example}

The canonical symplectic form (\ref{m83}) 
plays a fundamental role in view of
Darboux's theorem  \cite{libe}.

\begin{theorem}
Let $(Z,\Om)$ be a symplectic manifold. Every point $x$ of $Z$ has an
open neighborhood $U$ which is the domain of a coordinate chart 
$(y^1,\dots, y^n,p_1,\dots,p_n)$ such that the symplectic form $\Om$ has the
local expression (\ref{m83})
on $U$. Such coordinates are called {\it canonical}.
\end{theorem} 

\begin{proof} It is an immediate consequence of Proposition \ref{canpoiss}
and Proposition \ref{p11.5}.
\end{proof} 

\begin{remark}\label{canon}
Canonical coordinates on the manifold $T^*M$ are not adapted to the
fibration $T^*M\to M$ in general. For instance, the local coordinates
$({y'}^i=- p_i, p'_i=y^i)$ on $T^*M$ also are canonical.
\end{remark}

\newpage

\section{Polysymplectic geometry}

We consider first order Lagrangian and Hamiltonian formalisms on a
bundle $\Y$ over an $n$-dimensional base manifold $X$
\cite{cari,giac95,gunt,kijo,sard94,sard95}. 

\subsection{Lagrangian formalism}

Let $Y\to X$ be a bundle coordinatized by
$(x^\la, y^i)$. In jet terms, a {\it first order Lagrangian} is defined to
be a horizontal density $L=\cL\om$ 
on the jet manifold $J^1Y$ (see the notation (\ref{m112})). 
The jet manifold $J^1Y$ plays the role of the
finite-dimensional configuration  space of 
sections of $Y\to X$. We shall use the notation
\be
\pi^\la_i=\dr^\la_i\cL, \qquad \wt\pi=\cL-\pi^\la_iy^i_\la. 
\ee

We base our consideration on the {\it first variational 
formula} which provides the canonical decomposition of 
the Lie derivatives of
Lagrangians along projectable vector fields in accordance with the
variational task \cite{bau,fati,cam1,kru73,sard97}. 

Let $u=u^\la\dr_\la +u^i\dr_i$ be a projectable vector field on $\Y$
and $\ol u$ its jet lift (\ref{1.21}).
Given a Lagrangian density $L$, we have 
the following canonical decomposition
of the Lie derivative of $L$ along $u$:
\beq
\bx{\bL_uL\equiv
 u_V\rfloor \cE_L + h_0 d (u\rfloor\Xi_L),} \label{C30}
\eeq
where $u_V$ is the vertical part (\ref{31}) of $u$, 
\be
h_0: dy^i \mapsto y^i_\la dx^\la, \quad dy^i_\m\mapsto y^i_{\m\la} dx^\la,
\ee
is the operator of horizontalization,
\beq
\bx{\cE_L=(\dr_i- d_\la\dr^\la_i)\cL \,dy^i\w\om} \label{305}
\eeq
is the {\it Euler--Lagrange operator}, and $\Xi_L$ is some {\it Lepagean
equivalent} of $L$ on $J^1Y$. 

We restrict our consideration to the {\it Poincar\'e--Cartan form} 
\beq
\bx{\Xi_L=\pi^\la_idy^i\w\om_\la +\wt\pi\om.} \label{303}
\eeq
This is the only Lepagean equivalent which has the partner in
the framework of the Hamiltonian formalism (see Section 4.7). Moreover, 
if $n=1$, this is the unique Lepagean
equivalent of a Lagrangian.

The kernel $\Ker\cE_L\subset J^2Y$ of the Euler--Lagrange operator
(\ref{305})  
defines the system of second order {\it Euler--Lagrange equations}
\beq
\bx{(\dr_i- d_\la\dr^\la_i)\cL=0} \label{006}
\eeq
on the bundle $\Y$. On sections $s$ of $\Y$, these equations read
\beq
\dr_i\cL-(\dr_\la+\dr_\la s^j\dr_j
+\dr_\la\dr_\m s^j \dr^\m_j)\dr^\la_i\cL=0.\label{2.29}
\eeq

\begin{remark} Note that different Lagrangians  $L$ and $L'$ lead to the
same Euler--Lagrange operator iff
\beq
L'=L+h_0(\e), \label{m25}
\eeq
where $\e$ is a closed $n$-form on $Y$ \cite{kru73}.
Any closed form $\e$ on 
$Y$ is a Lepagean form. Let $L$ be a Lagrangian and
$\rho_L$ its Lepagian equivalent. Then, the Lepagian form $\rho_L+\e$ is the
Lepagean equivalent of the Lagrangian (\ref{m25}).
\end{remark}

\subsection{Legendre morphisms}

Every first order Lagrangian $L$ yields the {\it Legendre morphism}
$\wh L$ of the jet manifols $J^1Y$  to the {\it Legendre manifold}
\beq
\bx{\Pi=V^*Y\op\w_Y(\op\w^{n-1} T^*X)  =V^*Y\op\w_Y(\op\w^n
T^*X) \op\otimes_Y TX} \label{00}
\eeq
which plays the role of the finite-dimensional phase 
space of sections of $\Y$.
Given the bundle coordinates $(x^\la, y^i)$ of $\Y$,  
the Legendre bundle (\ref{00}) is coordinatized by  
$(x^\la ,y^i,p^\la_i)$, where $p^\la_i$ are the holonomic coordinates with
the transition functions
\beq
{p'}^\la_i = \det (\frac{\dr x^\ve}{\dr {x'}^\nu}) \frac{\dr y^j}{\dr{y'}^i}
\frac{\dr {x'}^\la}{\dr x^\m}p^\m_j.  \label{2.3}
\eeq
Relative to these coordinates, the Legendre morphism  $\wh L$ reads
\beq
 p^\m_i\circ\wh L=\pi^\m_i. \label{m11}
\eeq

The Poincar\'e--Cartan form $\Xi_L$ (\ref{303}) defines 
a morphism $\wh\Xi_L$ of the jet manifold $J^1Y$ to the {\it homogeneous
Legendre manifold}
\beq
\bx{Z = T^*Y\w(\op\w^{n-1}T^*X)}  \label{N41}
\eeq
provided with the holonomic coordinates $(x^\la, y^i, p^\la_i, p)$ with the
transition functions (\ref{2.3}) and
\beq
p'=\det (\frac{\dr x^\ve}{\dr {x'}^\nu})
(p-\frac{\dr y^j}{\dr {y'}^i}\frac{\dr {y'}^i}{\dr x^\mu}p^\mu_j).
\label{m10}
\eeq
Relative to these coordinates, the morphism $\wh\Xi_L$ reads
\beq
(p^\m_i, p)\circ\wh\Xi_L =(\pi^\m_i,\wt\pi ). \label{N42}
\eeq

 A glance at the expression (\ref{m10}) shows that $Z\to\Pi$ is a
1-dimensional affine bundle. We have the exact sequence 
\be
\bx{0\to \op\w^nT^*X\har Z\to\Pi \to 0.} 
\ee

The homogeneous Legendre manifold (\ref{N41}) is equipped with the canonical
$n$-form
\beq
\bx{\Xi= p\om + p^\la_i dy^i\w\om_\la.} \label{N43}
\eeq
Its coordinate expression (\ref{N43}) is 
maintained under holonomic coordinate
transformations (\ref{2.3}) and (\ref{m10}). The Poincar\'e--Cartan form
$\Xi_L$ (\ref{303}) is the pullback of $\Xi$ by the morphism $\wh\Xi_L$
(\ref{N42}).

\subsection{Polysymplectic structure}

The Legendre manifold (\ref{00})
possesses the {\it canonical polysymplectic form} 
 \beq
\bx{\bla =dp^\la_i\w dy^i\w\om\otimes\dr_\la} \label{406}
\eeq
whose coordinate expression (\ref{406}) is maintained under holonomic
coordinate transformations (\ref{2.3}). It is a pullback-valued form of
the type (\ref{1.11}).

\begin{remark} 
The  polysymplectic form (\ref{406}) can be introduced in
different ways. The Legendre manifold $\Pi$ is equipped also with 
the the {\it generalized Liouville form}
\beq
\bx{\bth =-p^\la_idy^i\w\om\otimes\dr_\la.} \label{2.4}
\eeq
Since (\ref{2.4}) is a pullback-valued form, one can not act on $\bth$
by the exterior differential in
order to recover the 
polysymplectic form $\bla$ (\ref{406}). At the same time, 
$\bla$ is the unique form which obeys the relation
\be
\bla\rfloor\f = -d(\bth\rfloor\f)
\ee
for any Pfaffian form $\f$ on $X$. 
\end{remark}

Given the atlas of holonomic coordinates $(x^\la, y^i, p^\la_i)$, let us
examine the coordinate transformations between these coordinates and any
coordinate atlas adapted to the bundle $\Pi\to X$ which keep invariant the
coordinate form (\ref{406}) of $\bla$. They will be called the {\it
polysymplectic canonical coordinate transformations}.

We find that,
since $y^i$ and $p^\la_i$ parameterize 
spaces of different dimensions if $n>1$,
polysymplectic canonical coordinate transformations  have a simpler
structure than that of the symplectic ($n=1$) ones (see Remark
\ref{canon}). Precisely they are compatible with the fibration
$\Pi\to Y$ and are exhausted by the
holonomic coordinate transformations (\ref{2.3}) and the translations 
\beq
{p'}^\la_i = p^\la_i + r^\la_i(y), \qquad \dr_jr^\la_i(y)=\dr_ir^\la_j(y).
\label{m242}
\eeq
Hereafter, we consider only holonomic coordinates 
$(x^\la ,y^i,p^\la_i)$ of $\Pi$.
 
\subsection{Hamiltonian connections}

Let $J^1\Pi$ be the jet manifold of the bundle
$\Pi\to X$. It is provided with the adapted  coordinates
$( x^\la, y^i, p^\la_i, y^i_\m, p^\la_{i\m})$.

\begin{definition}\label{d5.4} By analogy with the notion of a Hamiltonian
vector field (see Definition \ref{d11.2}), a connection
\be
\g =dx^\la\otimes(\dr_\la +\g^i_\la\dr_i +\g^\m_{i\la}\dr^i_\m)
\ee
on the bundle $\Pi\to X$ is said to be {\it locally
Hamiltonian}  [{\it Hamiltonian}] if the exterior form
$\g\rfloor\bla$ is closed [exact].
\end{definition}

It is readily observed that a connection $\g$ is locally Hamiltonian
iff it obeys the conditions
\beq
\dr^i_\la\g^j_\m-\dr^j_\m\g^i_\la=0,\quad
 \dr_i\g_{j\m}^\m- \dr_j\g_{i\m}^\m=0,\quad
 \dr_j\g_\la^i+\dr_\la^i\g_{j\m}^\m=0.\label{422}
\eeq

\begin{example}\label{e5.5} Given a linear symmetric connection $K$
(\ref{408}) on $T^*X$, every connection
$\G$ on the bundle $\Y$ gives rise to the connection
\beq
\wt\G =dx^\la\otimes[\dr_\la +\G^i_\la \dr_i +
(-\dr_j\G^i_\la  p^\m_i-K^\m{}_{\nu\la} p^\nu_j
 +K^\al{}_{\al\la} p^\m_j)\dr^j_\m]  \label{404}
\eeq
on $\Pi\to X$. It is easy to see that
$\wt\G\rfloor\bla =0$ and, consequently, the connection (\ref{404}) is a
locally Hamiltonian connection. Actually
$\wt\G$ appears to be a Hamiltonian connection (see Example \ref{HG}).
\end{example}

\subsection{Hamiltonian forms}

\begin{definition}\label{d5.8} A $n$-form $H$ on the
Legendre bundle $\Pi$ is called a 
{\it general Hamiltonian form} if
there exists a Hamiltonian connection such that
\be
\bx{\g\rfloor\bla = dH.} 
\ee
\end{definition}

Unless otherwise stated, general Hamiltonian forms 
will be considered modulo closed forms.

\begin{proposition}\label{p5.9} Let $H$ be a general Hamiltonian form. For
any horizontal density
$\wt H=\wt{\cH}\om$ on the bundle $\Pi\to X$, the form $H-\wt H$
is a Hamiltonian form.
\end{proposition}

The following example shows that general Hamiltonian forms on
$\Pi$ always exist. 

\begin{example}\label{HG}
Let $\G$ and $K$ be as in Example \ref{e5.5}. Then $\wt\G$ (\ref{404})	is a
Hamiltonian connection for the general Hamiltonian form
\be
\bx{H_\G =\G\rfloor\bth =p^\la_i dy^i\w\om_\la -p^\la_i\G^i_\la (y)\om,}
\ee
where $\bth$ is the generalized Liouville form (\ref{2.4}).
\end{example}

\begin{definition}
A general Hamiltonian form $H$ on $\Pi$ 
is said to be {\it Hamiltonian} if it 
has the splitting
\beq
\bx{H=H_\G -\wt H_\G=p^\la_idy^i\w\om_\la -
(p^\la_i\G^i_\la +\wt{\cH}_\G)\om
 =p^\la_idy^i\w\om_\la-\cH\om}  \label{4.7}
\eeq
modulo closed forms, where $\G$ is a
connection on $Y$ and $\wt H_\G$ is a horizontal density. 
\end{definition}

This splitting is preserved under the holonomic coordinate
transformations (\ref{2.3}), but not under translations (\ref{m242}). 

\begin{proposition}\label{p00}
There is the 1:1
correspondence between the Hamiltonian forms $H$ and the sections $h$ of the
bundle $Z\to \Pi$. We have 
\be
\bx{H=h^*\Xi,} 
\ee
where $\Xi$ is the canonical form (\ref{N43}) on $Z$
\end{proposition}

\begin{proof}
It is an immediate consequence of the expression 
(\ref{4.7}) for Hamiltonian forms. 
\end{proof}

By a {\it momentum morphism} we shall mean any bundle morphism 
\beq
\Phi:\Pi\op\to_Y J^1Y,\qquad \Phi= dx^\la\ot(\dr_\la +\Phi^i_\la\dr_i).
\label{m240}
\eeq
For instance, let $\G$ be a connection on the bundle $Y\to X$. Then,
the composition $\G\circ\pi_{\Pi Y}$ is a momentum morphism. Conversely,
every momentum morphism $\Phi$ defines
the associated connection $\G_\Phi =\Phi\circ\wh 0$
on $\Y$, where $\wh 0$ is the global zero section of $\Pi\to Y$.

\begin{proposition}\label{p5.12}
Every Hamiltonian form $H$ (\ref{4.7}) on the Legendre
manifold $\Pi$ yields the associated momentum morphism
\beq
 \wh H:\Pi\to J^1Y, \qquad (x^\la ,y^i,y_\la^i)\circ\wh H=( x^\la
,y^i,\dr^i_\la\cH), \label{415}
\eeq
and the associated connection $\G_H =\wh H\circ\wh 0$ on 
$Y\to X$. Conversely, every momentum morphism (\ref{m240}) defines the
Hamiltonian form
\be
H_\Phi =\Phi\rfloor\bth =p^\la_i dy^i\w\om_\la -p^\la_i\Phi^i_\la \om.
\ee
\end{proposition}

Given a Hamiltonian form $H$ (\ref{4.7}), there are the
algebraic conditions
\be
\g^i_\la =\dr^i_\la\cH, \qquad
\g^\la_{i\la}=-\dr_i\cH 
\ee
for a Hamiltonian connection $\g$ to be associated with a given
Hamiltonian form $H$. It should be emphasized that, if $n>1$, there 
exist different Hamiltonian connections for the same Hamiltonian form in
general. 

Let a Hamiltonian connection $\g$ associated with a Hamiltonian form $H$ have
an integral section $r$ of
$\Pi\to X$, that is, $\g\circ r=J^1r$. Then $r$ satisfies the
system of first order differential equations 
\bea
&& y^i_\la =\dr^i_\la\cH, \label{3.11a}\\
&&p^\la_{i\la} =-\dr_i\cH\label{3.11b}
\eea
on $\Pi$. They are called the {\it Hamilton equations}. 
It is readily observed that, 
if $r$ is a solution of the Hamilton equations (\ref{3.11a}) --
(\ref{3.11b}), it obeys the  relations
\be
J^1(\pi_{\Pi Y}\circ r)= \wh H\circ r. 
\ee

\subsection{Hamiltonian and Lagrangian formalisms}

We now turn to the relations between the Lagrangian formalism and the
polysymplectic Hamiltonian formalism. Let $L$ be a Lagrangian and $Q=\wh
L(J^1Y)$. We shall call $Q$ the {Lagrangian constraint space}. 

A  Hamiltonian form
$H$ is said to be {\it associated} with $L$ if it
obeys the conditions
\bea 
&&\wh L\circ\wh H|_Q=\Id_Q, \quad , 
\bx{p^\m_i=\dr^\m_i \cL(x^\la, y^j, 
\dr^j_\la\cH(p)), \quad p\in Q,}
\label{2.30a}\\
&& H_{\wh H}-H=L\circ\wh H, \qquad \bx{ p^\m_i\dr^i_\m\cH -\cH
\equiv\cL(x^\la, y^j, \dr^j_\la\cH).} \label{2.30b}
\eea

It should be emphasized that there may be different Hamiltonian forms
associated with $L$ in general.
 We restrict our consideration to Lagrangians which are {\it semiregular},
that is, the preimage $\wh L^{-1}(q)$ of each point 
$q\in Q$ is a connected submanifold of $J^1Y$.

 All Hamiltonian forms
associated with a semiregular Lagrangian $L$ coincide with each
other on the Lagrangian constraint space $Q$, and the Poincar\'e--Cartan form
$\Xi_L$ (\ref{303}) is the pullback 
\beq
\Xi_L=\wh L^*H,\qquad
\bx{\pi^\la_iy^i_\la-\cL\equiv \cH(x^\m,y^i,\pi^\la_i),} \label{Q2}
\eeq
of any such a Hamiltonian form $H$ by the Legendre morphism $\wh L$.
In this case, the following relation between solutions of the Euler--Lagrange
equations and solutions of the Hamilton equations takes place
\cite{sard93a,zakh}.
 
\begin{proposition}\label{pHL} 
Let a  section $r$ of $\Pi\to X$
be a solution of the Hamilton equations (\ref{3.11a}) -- (\ref{3.11b}) 
for a Hamiltonian form $H$ associated with a semiregular Lagrangian
$L$. If $r$ lives on the Lagrangian constraint space $Q$, the section
$ s=\pi_{\Pi Y}\circ r$ of $Y\to X$ satisfies the 
 Euler--Lagrange equations (\ref{2.29}) for $L$.
Conversely, given a semiregular Lagrangian $L$, let
$s$ be a solution of the
corresponding Euler--Lagrange equations.
Let $H$ be a Hamiltonian form associated with $L$ such that
\beq
\wh H\circ \wh L\circ J^1s=J^1s. \label{m210}
\eeq
Then, the section $r=\wh L\circ J^1s$ of $\Pi\to X$ is a solution of the
Hamilton equations (\ref{3.11a}) -- (\ref{3.11b})  for $H$.
\end{proposition}
 
We say that a family of Hamiltonian forms $H$
associated with a semiregular Lagrangian $L$ is
{\it complete} if, for each solution $s$ of the Euler--Lagrange
equations, there exists a
solution $r$ of the Hamilton equations  for
some  Hamiltonian form $H$ of this family so that
\beq
r=\wh L\circ J^1s, \qquad
s= \pi_{\Pi Y}\circ r. \label{2.37}
\eeq

\begin{example}\label{regular}
In case of a {\it hyperregular} Lagrangian $L$ (i.e., the Legendre
morphism
$\wh L$ is a diffeomorphism), the Lagrangian formalism and the
polysymplectic Hamiltonian formalism are equivalent.
There exists the unique Hamiltonian form
\be
H=H_{\wh L^{-1}}+L\circ\wh L^{-1}
\ee
associated with $L$. The corresponding momentum morphism (\ref{415}) is the
diffeomorphism
$\wh H=\wh L^{-1}$. As a consequence, there is the 1:1 correspondence
between the solutions of the Euler--Lagrange equations for $L$ 
and the Hamilton equations for $H$.
In case of a {\it regular} Lagrangian $L$ (i.e., $\wh L$ is a local
diffeimorphism), the Lagrangian constraint space $Q$
is an open submanifold of the Legendre manifold $\Pi$. 
If a regular Lagrangian density is also
semiregular, the associated Legendre
morphism is a diffeomorphism of $J^1Y$ onto $Q$ and, on $Q$, we	can recover
all results true for  hyperregular Lagrangians.
\end{example}

\begin{remark}\label{joint}
Given a Hamiltonian form $H$ (\ref{4.7}) on the Legendre manifold $\Pi$
(\ref{00}), let us consider the Lagrangian 
\beq
\bx{L_H = (p^\la_iy^i_\la - \cH)\om} \label{Q3}
\eeq
on the configuration space $J^1\Pi$ coordinatized by $(x^\la,y^i, p^\m_i,
y^i_\la, p^\m_{i\la})$. This Lagrangian does not depend on the
coordinates $p^\m_{i\la}$.
It is readily observed that
the Poincar\'e--Cartan form $\Xi_{L_H}$ (\ref{303})
of the Lagrangian (\ref{Q3}) consists with the Hamiltonian form $H$ and
the Euler--Lagrange equations (\ref{006}) for $L_H$ recover the Hamilton
equations (\ref{3.11a}) -- (\ref{3.11b}) for $H$. 
\end{remark}

\subsection{Vertical extension of the polysymplectic formalism}

In time-dependent mechanics, the
machinery that we present below provides the the way to maintain the form 
(\ref{4.7}) of Hamiltonian forms under canonical
transformations. By analogy with the BRS  generalization of mechanics
\cite{gozz89,gozz92}, it represents also a first step toward the BRS
quantization of the polysymplectic Hamiltonian formalism. 

Given a bundle $\Y$, let us consider its vertical tangent bundle $VY$
coordinatized by $(x^\la, y^i,\dot y^i)$. We show that the Hamiltonian
formalism for sections of $\Y$ is naturally extended to the Hamiltonian
formalism for sections of $VY\to X$.

The Legendre manifols (\ref{00}) corresponding to $VY\to X$ is
\be
\Pi_{VY}=V^*VY\op\w_{VY}(\op\w^{n-1} T^*X). 
\ee
It is coordinatized by $(x^\la, y^i, \dot y^i, q^\la_i, v^\la_i)$.

\begin{proposition}\label{vertmom}
In virtue of the bundle isomorphism (\ref{m6}), there exists the
bundle isomorphism
\beq
\bx{\Pi_{VY}\op=_{VY}V\Pi,} 
\qquad q^\la_i\llra\dot p^\la_i, \qquad v^\la_i\llra p^\la_i, \label{m14}
\eeq
where $(x^\la, y^i,  p^\la_i, \dot y^i, \dot p^\la_i)$ are the coordinates of 
$V\Pi$.
\end{proposition}

We shall utilize the compact notation
\beq
\dot\dr_i=\frac{\dr}{\dr\dot y^i}, \qquad \dot\dr^i_\la=\frac{\dr}{\dr\dot
p_i^\la}. \label{m147}
\eeq
 Recall also the notation $\dr_V$ (\ref{m141}).

One can develop the
Hamiltonian formalism on $\Pi_{VY}$ by analogy with that on $\Pi$. The
manifold
$V\Pi$ is equipped with the canonical polysymplectic form 
\beq
\bx{\bla_V=[d\dot p^\la_i\w dy^i +dp^\la_i\w d\dot y^i]\w\om\ot\dr_\la.}
\label{m16}
\eeq
Its coordinate expression is maintained under 
holonomic transformations of the
composite bundle $V\Pi\to\Pi\to Y$.

\begin{proposition}\label{p01}
Let $\g$ be a Hamiltonian connection on $\Pi$ associated with a Hamiltonian
form $H$ (\ref{4.7}).
Then, the vertical connection $V\g$ (\ref{43}) is a Hamiltonian connection
associated with the Hamiltonian form
\beq
\bx{H_V=(\dot p^\la_idy^i -\dot y^idp^\la_i)\w\om_\la -\cH_V,
\quad \cH_V =\dr_V\cH=(\dot y^i\dr_i +\dot p^\la_i\dr^i_\la)\cH.} 
\label{m17}
\eeq
\end{proposition}

\begin{proof} It is easily justified that, given a Hamiltonian connection
\be
\g=dx^\m\ot (\dr_\m +\g^i_\m\dr_i +\g^\la_{i\m}\dr^i_\la),
\quad
\g^i_\m=\dr^i_\m\cH, \qquad \g^\la_{i\la} =-\dr_i\cH,
\ee
the vertical connection
\be
V\g=dx^\m\ot [\dr_\m +\g^i_\m\dr_i +\g^\la_{i\m}\dr^i_\la
+ \dr_V\g^i_\m\dot\dr_i +\dr_V\g^\la_{i\m}\dot\dr_\la^i]
\ee
obeys the Hamilton equations for the Hamiltonian form (\ref{m17}):
\be
&& \g^i_\m=\dot\dr_\m^i\cH_V =\dr^i_\m\cH,\\
&& \g^\la_{i\la}=-\dot\dr_i\cH_V =-\dr_i\cH,\\
&& \dot \g^i_\m=\dr^i_\m\cH_V=\dr_V\dr^i_\m\cH,\\
&& \dot \g^\la_{i\la}=-\dr_i\cH_V=-\dr_V\dr_i\cH.
\ee
\end{proof}

In particular, if there is the splitting $\cH=p^\la_i\G^i_\la +\wt\cH$
relative to some connection $\G$ on $\Y$, then we have the splitting
\be
\cH_V=\dot p^\la_i\G^i_\la -\dot y^j(-p^\la_i\dr_j\G^i_\la) +\dr_V\wt\cH
\ee
with respect to the lift $\wt\G$ (\ref{404}) of $\G$ onto $\Pi\to X$.

Note that the Hamiltonian form $H_V$ (\ref{m17}) 
can be obtained also in the
following way. Given the homogeneous Legendre manifold $Z$ (\ref{N41}), let
us consider the vertical tangent bundle $VZ$ of $Z\to X$ coordinatized by
$(x^\la, y^i,  p^\la_i, p, \dot y^i, \dot p^\la_i, \dot p)$. 
It is provided with the canonical form
\be
\bx{\Xi_V=\dot p\om + \dot p^\la_i dy^i\w\om_\la - \dot y^i dp^\la_i
\w\om_\la} 
\ee
whose expression is maintained under holonomic coordinate transformations.
Note that one can utilize also the form $\Xi_V + 
d(\dot y^i p^\la_i) \w\om_\la$
since the form $d(\dot y^i p^\la_i) \w\om_\la$ is well-behaved.

Put $H=h^*\Xi$, where $h$ is a section of the bundle $Z\to \Pi$. 
Then, we have 
\be
\bx{H_V=(Vh)^*\Xi_V,} 
\ee
where $Vh:V\Pi\to VZ$ is the vertical tangent morphism to $h$.

We now turn to the vertical extension of the Lagrangian formalism on
$J^1Y$ onto the
configuration space $VJ^1Y=J^1VY$ coordinatized by $(x^\la, y^i, y^i_\la,
\dot y^i, \dot y^i_\la)$. Given a Lagrangian $L$ on $J^1Y$, 
let us consider the Lagrangian
\beq
 L_V=\pr_2\circ VL: VJ^1Y\to \op\w^nT^*X, \quad
 \bx{\cL_V= \dr_V\cL=(\dot y^i\dr_i + \dot y^i_\la\dr_i^\la)\cL,} 
\label{m18}
\eeq
on $VJ^1Y$. It is readily observed that the variational derivative
$\dot\dl_i \cL_V=\dl_i\cL$
recovers the Euler--Lagrange equations (\ref{006}).

The Lagrangian (\ref{m18}) yields the Legendre morphism
\ben
&& \wh L_V=V\wh L: VJ^1Y\op\to_{VY} V\Pi, \label{m19}\\
&& p^\la_i=\dot\dr^\la_i\cL_V=\pi^\la_i, \qquad \dot p^\la_i =
\dr_V\pi^\la_i. \nonumber
\een
Conversely, given a Hamiltonian form $H_V$ 
(\ref{m17}) on $V\Pi$, there is the
momentum morphism
\be
&& \wh H_V=V\wh H: V\Pi\op\to_{VY} VJ^1Y, \\
&& y^i_\la=\dot\dr^i_\la\cH_V =\dr^i_\la\cH, \qquad
\dot y^i_\la= \dr_V\dr^i_\la\cH. 
\ee

Let us consider the relation between the Hamiltonian form
$H_V$ and the Lagrangian $L_V$ if the 
Hamiltonian form $H$ is associated with
the Lagrangian $L$.

\begin{proposition}\label{p02} The Legendre morphism (\ref{m19}) is a
surjection of $VJ^1Y$ onto $VQ$. 
\end{proposition}

\begin{proof} One can show that the equations 
\be
\dot p^\la_i =
(\dot y^i\dr_i + \dot y^i_\la\dr_i^\la)\pi^\la_i
\ee
are equivalent to the equations  
\be
(\dot p^\la_i\dr_\la^i + \dot y^i\dr_i)\rfloor d[p^\m_i-\dr^\m_i \cL(x^\la,
y^j, \dr^j_\la\cH(p))]=0
\ee
characterizing tangent vectors to the fibres of the Lagrangian constraint
bundle $Q$. 
\end{proof}

Moreover, $VQ$ appears to be the image of $\wh L_V$ restricted to 
$\wh H (Q)$. It follows that a relation similar to (\ref{2.30a})
takes place. At the same time, a relation similar to
(\ref{2.30b}) holds only on the constraint space $Q$. 
 Let a Hamiltonian form $H$ be associated with a
semiregular Lagrangian $L$. Then, the Hamiltonian form $H_V$ and the
Lagrangian $L_V$ (which fails to be semiregular in general) satisfy 
the relation similar to (\ref{Q2}) on $\wh H(Q)$.

\newpage

\section{Time-dependent Hamiltonian mechanics}

To describe time-dependent mechanical systems, let us consider a bundle $\Y$
with a $m$-dimensional standard fibre $M$ over 
a 1-dimensional base $X$. It is
provided with bundle coordinates $(t,y^i)$. Observe that: 

(i) the jet manifold $J^1Y$ is modelled on the 
vertical tangent bundle $VY$ of
$Y$;

(ii) the Legendre
bundle $\Pi$ (\ref{00}) over $Y$ is the vertical cotangent bundle $V^*Y$ of
$Y$ coordinatized by $(t,y^i,p_i)$;

(iii) the homogeneous Legendre bundle $Z$ (\ref{N41}) over $Y$ is 
the cotangent bundle
$T^*Y$ of $Y$ coordinatized by $(t,y^i,p,p_i)$.

\begin{remark}
If the base manifold is
contractible, i.e. $X= \R$, the bundle $\Y$ is trivializable. Given a
trivialization
\beq
Y\simeq\R\xx M, \label{m33}
\eeq
we have the corresponding splittings
\ben
&&J^1Y\simeq\R\xx TM \nonumber\\
&&\Pi\simeq\R\xx T^*M. \label{m35}
\een
\end{remark}

\subsection{$n=1$ Reduction of the polysymplectic formalism}

The  phase space $\Pi=V^*Y$. It is provided with
the holonomic coordinates $(t,y^i,p_i)$
possessing the transition functions 
\beq
p'_i=\frac{\dr y^j}{\dr {y'}^i}p_j. \label{m153}
\eeq

The Legendre manifold $\Pi$ admits
the canonical form $\bla$ (\ref{406}) which reads
\beq
\bx{\bla=dp_i\w dy^i\w dt\ot\dr_t.} \label{m36}
\eeq
As a particular case of the polysymplectic machinery of the previous Section,
we say that a connection
\be
\g =dt\ot(\dr_t +\g^i\dr_i +\g_i\dr^i)
\ee
on the bundle $\Pi\to X$ is locally 
Hamiltonian [Hamiltonian] if the exterior form
$\g\rfloor\bla$ is closed [exact]. A connection $\g$ is 
locally Hamiltonian iff it obeys the conditions (\ref{422}) which
now take the form
\be
\dr^i\g^j-\dr^j\g^i=0,\quad
\dr_i\g_j- \dr_j\g_i=0,\quad \dr_j\g^i+\dr^i\g_j=0.
\ee

As in Example \ref{e5.5}, we observe that every connection
$\G=dt\ot(\dr_t +\G^i\dr_i)$ on the bundle
$\Y$ gives rise to the Hamiltonian connection
\beq
\bx{\wt\G =dt\ot(\dr_t +\G^i\dr_i -\dr_j\G^i p_i\dr^j)}  \label{m38}
\eeq
on $\Pi\to X$ which consists with the covertical connection $V^*\G$
(\ref{44}). The corresponding Hamiltonian form is
\beq
H_\G=p_idy^i -p_i\G^idt. \label{m61}
\eeq

Let $H$ be a Hamiltonian form (\ref{4.7}) on $\Pi=V^*Y$. It reads
\beq
\bx{H=p_idy^i-\cH dt=p_idy^i -p_i\G^idt -\wt{\cH}_\G dt.}  \label{m46}
\eeq
We call $\cH$ and $\wt\cH$ in the decomposition (\ref{m46}) the
{\it Hamiltonian} and the {\it Hamilton function} respectively.
Let $\g$ be a Hamiltonian connection on $\Pi\to X$ associated with the
Hamiltonian form (\ref{m46}). It satisfies the relations
\ben
&&\g\rfloor\bla =dp_i\w dy^i+ \g_idy^i\w dt -\g^idp_i\w dt = dH,\nonumber\\
&&\g^i =\dr^i\cH, \qquad \g_i=-\dr_i\cH. \label{m40}
\een
A glance at the equations (\ref{m40}) shows that, in the case of mechanics,
one and only one Hamiltonian connection is associated with a given 
Hamiltonian form. 

In accordance with Remark \ref{coneq}, every connection $\g$ on $\Pi\to X$
yields the system of first order 
differential equations (\ref{39}) which reads
\beq
y^i_t =\g^i, \qquad p_{it} =\g_i.\label{m170}
\eeq
They are called the {\it evolution equations}. 
If $\g$ is a Hamiltonian connection associated with the Hamiltonian form $H$
(\ref{m46}), the evolution equations (\ref{m170}) come to the
Hamilton equations 
\bea
&& y^i_t =\dr^i\cH, \label{m41a}\\
&& p_{it} =-\dr_i\cH.\label{m41b}
\eea

Note that, once a
trivialization (\ref{m35}) is chosen, the Hamiltonian form (\ref{m46}) 
yields the well-known Poincar\'e--Cartan integral invariant \cite{libe}.
At the same time, the splitting (\ref{m46}) is not
maintained under canonical 
transformations (see Section 5.4). This fact calls
into play the general Hamiltonian forms (see Proposition \ref{genform}).

Another well-known ingredient in time-dependent mechanics is 
the horizontal lift 
\beq
\bx{\tau_H=\dr_t +\dr^i\cH\dr_i -\dr_i\cH\dr^i} \label{m57}
\eeq
onto $\Pi$ of the standard vector field $\dr_t$ on $X$ by means of a
Hamiltonian connection $\g$ associated with a Hamiltonian form $H$
(\ref{m46}).
It is a nowhere vanishing vector field on $\Pi$ which
obeys the relations
\beq
\tau_H\rfloor H =p_i\dr^i\cH-\cH, \qquad \tau_H\rfloor dH=0. \label{m58}
\eeq
We call $\tau_H$ (\ref{m57}) the {\it horizontal Hamiltonian vector field} of
the Hamiltonian form $H$. 

\begin{remark}\label{r13.1}
Every connection $\g$ on a bundle $\Pi\to X$ is a curvature-free connection
(see Example \ref{rcon}).
In virtue of Proposition \ref{flat}, such a connection defines a horizontal
foliation on
$\Pi\to X$. Its leaves are the integral curves of the horizontal lift
\beq
\tau_\g=\dr_t + \g^i\dr_i +\g_i\dr^i \label{m115}
\eeq
of $\dr_t$ by $\g$. The corresponding Pfaffian
system is locally generated by the forms  $(dy^i-\g^i dt)$
and $(dp_i-\g_i dt)$.
\end{remark}

It follows that every Hamiltonian
connection and, accordingly, every Hamiltonian form defines the 
corresponding Hamiltonian foliation on $\Pi$.
Its leaves are integral curves of the horizontal Hamiltonian vector field
(\ref{m57}). One can think of these integral curves as being the generalized
solutions of the Hamilton equations (\ref{m41a}) and (\ref{m41b}) (in
accordance with the definition of generalized solutions given in
\cite{kras}). They locally coincide with the integral sections of the
Hamiltonian connection $\g$.

Given a function $f$ on
$\Pi$, we have the {\it Hamiltonian evolution equation} 
\beq
\bx{\tau_H\rfloor df=d_{Ht}f=(\dr_t +\dr^i\cH\dr_i -\dr_i\cH\dr^i)f}
\label{m59}
\eeq
relative to the Hamiltonian $\cH$. On solutions $r$ of the Hamilton
equations, (\ref{m59}) is equal to the total time
derivative of the function $f$:
\be
r^*d_{Ht}f=\frac{d}{dt}(f\circ r).
\ee

The goal is to write the Hamiltonian evolution equation (\ref{m59}) in the
terms of a Poisson bracket.

\subsection{Canonical Poisson structure}

Let us consider the homogeneous  phase space $Z=T^*Y$. The canonical form
$\Xi$ (\ref{N43}) on it comes to the canonical Liouville form
\beq
\bx{\Xi=pdt +p_idy^i.} \label{m91}
\eeq
Its exterior differential is the canonical symplectic form
\beq
\bx{d\Xi=dp\w dt +dp_i\w dy^i.} \label{m92}
\eeq
The corresponding Poisson bracket on the space $S(Z)$ of
functions on $Z$ reads
\beq
\bx{\{f,g\} =\dr^pf\dr_tg - \dr^pg\dr_tf +\dr^if\dr_ig-\dr^ig\dr_if.}
\label{m116}
\eeq

Let us consider the subspace of $S(Z)$ which
consists of the pullbacks of functions on $\Pi$ by the projection $Z\to
\Pi$. It is easily observed that this subspace is closed under the Poisson
bracket (\ref{m116}). Then,
according to Proposition \ref{p11.3}, one can 
show that the canonical Poisson
structure (\ref{m116}) on $Z$  induces the canonical Poisson structure 
\beq
\bx{\{f,g\}_V =\dr^if\dr_ig-\dr^ig\dr_if}
\label{m72}
\eeq
on $\Pi$ by the projection $Z\to \Pi$. The corresponding bivector on 
$\Pi$ is vertical with respect to the projection $\Pi\to X$. It reads
\be
w^{ij} =0, \qquad w_{ij}=0, \qquad w^i{}_j=1. 
\ee
 Since the rank of $w$ is constant, the Poisson structure (\ref{m72}) is
regular.

The Poisson structure (\ref{m72}) is obviously degenerate. It defines the
symplectic foliation on $\Pi$ which coincides with the fibration
$\Pi\to X$. The Hamiltonian 
vector fields associated with the Poisson bracket (\ref{m72}) are the
vertical vector fields on $\Pi\to X$. The Hamiltonian vector field $\vt_f$ 
of a function
$f$ is defined by the relation (\ref{m80}):
\be
\{f,g\}_V=\vt_f\rfloor dg, \qquad g\in S(\Pi).
\ee
It reads
\beq
\bx{\vt_f = \dr^if\dr_i- \dr_if\dr^i.} \label{m73}
\eeq

Note that the bundle coordinates $(t,y^i,p_i)$ of $\Pi$ are exactly the
canonical coordinates (\ref{m111}) for the Poisson structure (\ref{m72}). In
particular, the symplectic forms on 
the fibres of $\Pi\to X$ are the pullbacks 
\be
\Om_x=dp_i\w dy^i
\ee
of the canonical symplectic form on the standard fibre
$T^*M$ of the bundle $\Pi\to X$ with respect to morphisms of trivialization.

The Poisson structure (\ref{m72}) on $\Pi$ can be introduced in a different
way \cite{cari89}. There exists the canonical closed 3-form
\beq
\bx{\bom=dp^i\w dy_i\w dt,} \label{m66}
\eeq
on the Legendre manifold $\Pi$. 
With this form, every function
$f$ on $\Pi$ defines a vertical vector field $\vt_f$ 
on the bundle $\Pi\to X$
by the relation
\be
\bx{\vt_f\rfloor\bom = df\w dt.} 
\ee
Then, the Poisson bracket (\ref{m72}) is given by condition
\beq
\vt_g\rfloor\vt_f\rfloor\bom=\{f,g\}_Vdt. \label{m124}
\eeq

The canonical forms $\bla$ (\ref{m36}) and $\bom$ (\ref{m66}) on
$\Pi$ can be considered on the same footing as follows.

\begin{proposition}\label{p13.2} Let $u$ be a vector field on $\Pi$
projected onto the standard vector field $\dr_t$ on $X$. This vector
field obeys the relation
\beq
\bL_u\bom=d(u\rfloor\bom)=0 \label{m68}
\eeq
iff it is the horizontal lift $\tau_\g$ (\ref{m115})
of $\dr_t$ onto $\Pi$ by means of a locally Hamiltonian connection $\g$ on
$\Pi\to X$. In particular, $\tau_\g\rfloor\bom=dH$ 
if $\g$ is a Hamiltonian connection and $H$ is the corresponding
general Hamiltonian form.
\end{proposition}

\begin{proposition}\label{p13.4} 
If $\g$ is a Hamiltonian connection associated
with a Hamiltonian form $H$ and $\vt_f$ is the Hamiltonian vector
field (\ref{m73}), then $\g +\vt_fdt$ is a Hamiltonian connection
associated with the Hamiltonian form $H +fdt$.
\end{proposition}

Given a Hamiltonian connection $\g$ for a Hamiltonian form $H$, let us
consider its splitting 
\beq
H=H_0+\wt\cH dt, \qquad \g=\g_0 +\vt dt, \label{m95}
\eeq
where $H_0$ is some Hamiltonian form, 
$\g_0$ is the Hamiltonian connection for $H_0$, and $\vt$ is the Hamiltonian 
vector field of the function $\wt\cH$. 
One can bring the Hamiltonian evolution equation (\ref{m59}) relative to
$H$ into the form compatible with the splitting (\ref{m95}). It reads
\beq
d_{Ht}f= d_{H_0t}f +\{\wt\cH,f\}_V= \dr_tf +(\g_0{}^i\dr_i +\g_{0i}\dr^i)f
+\{\wt\cH,f\}_V. \label{m96}
\eeq
A glance at this expression shows that Hamiltonian evolution equations in
time-dependent mechanics do not reduce to the Poisson brackets. 

\begin{imp} This fact
becomes relevant to the quantization problem. The second term on the
right-hand  side of the equation (\ref{m96}) remains classical.\end{imp}

In this context, the main problem is to express 
the Hamiltonian evolution equation of a classical 
system in terms of the Poisson bracket. Then, one can
bring this  Hamiltonian evolution equation into the operator evolution 
equation under quantization. A glance at the expression (\ref{m96}) 
shows that this is possible only with respect to the splitting (\ref{m95}),
where the connection
$\g_0$ is brought into zero by a canonical coordinate transformation (see
Section 5.4).

{} From physical viewpoint, the splitting (\ref{m95}) has a meaning if the
connection $\g_0$ characterizes a reference frame (see Section 5.6). 

\subsection{Presymplectic and contact structures}

Besides the canonical Poisson structure, there is no other canonical
structure on the  phase space $\Pi=V^*Y$ of time-dependent mechanics in
general. At the same time, there are structures on $\Pi$ specified by the
choice of a Hamiltonian form $H$.

In virtue of Proposition \ref{p00}, every Hamiltonian form
$H$ on the  phase space $\Pi$ is the pullback $H=h^*\Xi$ of the Liouville
form (\ref{m91}) by a section $h$ of the bundle $Z\to \Pi$.
Accordingly, its differential 
\be
dH=(dp_i+\dr_i\cH dt)\w(dy^i-\dr^i\cH dt) 
\ee
is the pullback $h^*d\Xi$ of the symplectic
form (\ref{m92}). 
It is a presymplectic form of the constant rank
$2m$ since the form
\beq
(dH)^m =(dp_i\w dy^i)^m -m(dp_i\w dy^i)^{m-1}\w d\cH\w dt \label{m90}
\eeq
is obviously nowhere vanishing. However, this presymplectic structure does
not introduce any essentially new object because the corresponding 
Hamiltonian
vector fields are proportional to the horizontal Hamiltonian vector field
$\tau_H$ (\ref{m57}).
At the same time, a Hamiltonian form (\ref{m46}) satisfying certain
conditions is a contact form which defines a nondegenerate Jacobi structure
on $\Pi$ as follows.

\begin{proposition}\label{p13.1} A Hamiltonian form (\ref{m46}) is a contact
form if the density 
\be
[\cH]=p_i\dr^i\cH-\cH 
\ee
nowhere vanishes \cite{libe}.
\end{proposition}

\begin{proof} Since the horizontal Hamiltonian vector field $\tau_H$
(\ref{m58}) is nowhere vanishing, the condition $H\w(dH)^m\neq 0$ is
equivalent to the condition 
\be
\tau_H\rfloor(H\w(dH)^m)= (\tau_H\rfloor H) (dH)^m =[\cH]dH)^m\neq 0
\ee
and, since the form $(dH)^m$ (\ref{m90}) is nowhere vanishing, the result
follows.
\end{proof}

\begin{remark}
In order to make $[\cH]$ nowhere vanishing, one may add some exact
form (e.g., $cdt$, $c=$const.) to $H$. For instance, the Hamiltonian form
$H_\G$ (\ref{m61}) is not a contact form because $[\cH]=0$, but the
equivalent form $H_\G- dt$, where $[\cH]=1$, is it.
\end{remark} 

Given a Hamiltonian form $H$, let $[\cH]$ be nowhere vanishing so that 
the form $H$ is a contact form. The corresponding Reeb vector field reads 
\beq
E_H= [\cH]^{-1}\tau_H. \label{m62}
\eeq
In virtue of Proposition
\ref{p11.2}, this form has the associated Jacobi bracket defined by the Reeb
vector field (\ref{m62}) and by the bivector field $w_H$ derived from the
relations (\ref{m55}). We find
\be
&& w_H(\f,\si)= w_H^\sh\f\rfloor\si=\f^i\si_i + p_i\si^iE_H\rfloor\f - [\f
\llra\si], \\
&& w_H^\sh\f=-p_i\f^i [\cH]^{-1}\dr_t +(\f^i- p_j\f^j
[\cH]^{-1}\dr^i\cH)\dr_i + \\
&&\qquad (-\f_i + [\cH]^{-1}(p_j\f^j\dr_i\cH+
p_i\tau_H\rfloor\f))\dr^i,
\ee
where $\f$ and $\si$ are Pfaffian forms on $\Pi$.
The corresponding Jacobi bracket (\ref{m50}) reads
\beq
\{f,g\}_H= \{f,g\}_V + [\cH]^{-1}([g]d_{Ht}f - [f]d_{Ht}g),
\label{m64}
\eeq
where $\{f,g\}_V$ is the canonical Poisson bracket (\ref{m72}) and 
\be
 [f]=p_i\dr^if-f, \qquad [g]=p_i\dr^ig-g.
\ee
In particular, let $H$ have the splitting (\ref{m95}). We find
\beq
\{\wt\cH,g\}_H =  [\cH]^{-1}([g]d_{H_0t}\wt\cH - [\wt\cH]d_{H_0t}g).
\label{m130}
\eeq

Given a contact Hamiltonian form $H$, one can consider also the Jacobi
bracket defined by the Reeb vector field $E_H$ (\ref{m62}) alone. It reads
\beq
\{f,g\}_E= [\cH]^{-1}(fd_{Ht}g -gd_{Ht}f) \label{m132}.
\eeq
Given the splitting (\ref{m95}), we find
\beq
\{\wt\cH,g\}_E =  [\cH]^{-1}(\wt\cH d_{Ht}g - gd_{H_0t}\wt\cH). 
\label{m131}
\eeq

A glance at the expressions (\ref{m130}) and (\ref{m131}) shows that the
Jacobi brackets (\ref{m64}) and (\ref{m132}) have no advantage over the
Poisson bracket (\ref{m72}) in order to write the Hamiltonian evolution
equation (\ref{m59}).

\subsection{Canonical transformations}

Up to now, we have followed the 
polysymplectic Hamiltonian formalism and have 
considered transformations which keep the fibration $\Pi\to Y$. Let us
now examine canonical transformations of 
time-dependent mechanics which are not
compatible with this fibration. Remind that the base $X$ is
not transformed.

\begin{definition}\label{can1} Given an atlas $\Psi=\{\psi_\xi\}$ of the
bundle $\Pi\to X$, the bundle coordinates $(t,y^i,p_i)$, where
\be
y^i(p)=(y^i\circ\pr_2\circ \psi_\xi)(p), \quad p_i(p)=(p_i\circ\pr_2\circ
\psi_\xi)(p), \quad p\in \Pi,
\ee
are said to be the {\it canonical coordinates} if, in these
coordinates, the form
$\bla$ and equivalently the form $\bom$  are given by the
canonical expressions (\ref{m36}) and (\ref{m66}).
\end{definition}

The {\it canonical coordinate
transformations} satisfy the relations
\ben
&& \frac{\dr {p'}_i}{\dr p_j}\frac{\dr {y'}^i}{\dr p_k}
-\frac{\dr {p'}_i}{\dr p_k}\frac{\dr {y'}^i}{\dr p_j}=0, \nonumber\\
&& \frac{\dr {p'}_i}{\dr y^j}\frac{\dr {y'}^i}{\dr y^k}
-\frac{\dr {p'}_i}{\dr y^k}\frac{\dr {y'}^i}{\dr y^j}=0, \label{m44}\\
&& \frac{\dr {p'}_i}{\dr p_j}\frac{\dr {y'}^i}{\dr y^k}
-\frac{\dr {p'}_i}{\dr y^j}\frac{\dr {y'}^i}{\dr p_k} = \dl^k_j. \nonumber
\een

By definition, the
holonomic coordinates of $\Pi=V^*Y$ are obviously canonical coordinates.

\begin{definition}\label{can2} By a
{\it canonical transformation (morphism)} is meant an isomorphism $\rho$ of
the bundle
$\Pi\to X$ over $X$ such that any atlas $\Psi$ of holonomic coordinates
$(t,y^i,p_i)$ of $\Pi$ and the atlas
$\Psi\circ\rho^{-1}$ (\ref{m136}) of the coordinates
\be
(t,{y'}^i=y^i\circ \rho^{-1},{p'}_i=p_i\circ \rho^{-1}) 
\ee
are related by the canonical coordinate transformations (\ref{m44}).
\end{definition}

The equivalent coordinate-free definition of a canonical morphism is the
following.

\begin{definition}\label{can3} A
canonical morphism is an isomorphism $\rho$ of
the bundle $\Pi\to X$ over $X$
which preserve the canonical form $\bom$ (\ref{m66}), that is,
$\rho^*\bom=\bom$.
\end{definition}

It is easily observed that canonical morphisms preserve the
canonical Poisson structure (\ref{m72}) on $\Pi$, that is, 
\be
\{f\circ\rho,g\circ\rho\}_V=(\{f,g\}_V)\circ\rho.
\ee

\begin{proposition}\label{genham}
Canonical morphisms send Hamiltonian connections to Hamiltonian connections. 
\end{proposition}

\begin{proof}
The proof is based on the relation $T\rho(\tau_\g) =\tau_{\rho(\g)}$,
where
$\g$ is a connection on $\Pi\to X$ and $\tau_\g$ is the horizontal vector
field (\ref{m115}). If $\g$ is a Hamiltonian connection such that
$\tau_\g\rfloor\bom =dH$, we have 
\be
\tau_{\rho(\g)}\rfloor\bom=(\rho^{-1})^*(\tau_\g\rfloor\bom)
=d((\rho^{-1})^*H).
\ee
\end{proof}

A glance at the relation (\ref{m68}) shows that, for each locally
Hamiltonian connection $\g$, the horizontal 
Hamiltonian vector field $\tau_\g$
is the generator of a local 1-parameter group 
$G_\g$ of canonical morphisms of
$\Pi\to X$. It leads to the following assertion.

\begin{proposition}\label{p13.5} Let $X=\R$ and $\g$ be a complete locally
Hamiltonian connection on $\Pi\to \R$. There exist canonical coordinate
transformations which bring $\g$ into zero.
\end{proposition}

\begin{proof}
In virtue of Proposition \ref{complcon}, 
there exists a trivialization of the
bundle $\Pi\to \R$ such that $\g^i=0$, $\g_i=0$ relative to 
coordinates which are 
constant along the integral curves of $\tau_\g$. Since $G_\g$ is a group
of canonical transformations, we deduce that the above-mentioned coordinates
are canonical.
\end{proof}

\begin{imp} From physical viewpoint, the above coordinates are the
initial values of the canonical variables.\end{imp} 

\begin{corollary}\label{c13.5}
The evolution equations (\ref{m170}) 
associated with a Hamiltonian connection
can be locally brought into the equilibrium equations
\be
y^i_t=0, \qquad p_{it}=0
\ee
by canonical transformations. 
\end{corollary}

\begin{example}\label{prim1}
Let us consider 1-dimensional motion with constant accellaration $a$ with
respect to a reference frame whose coordinates are $(t,y)$. The
corresponding Hamiltonian and the Hamiltonian connection read
\ben
&&\cH=\frac{p^2}{2} -ay, \nonumber \\
&& \g^y= p, \qquad \g_p =a. \label{m157}
\een
This is a complete connection. The canonical transformation
\be
y'=y-pt +\frac{at^2}{2}, \qquad p'=p-at 
\ee
brings the connection (\ref{m157}) into zero.
\end{example}

\begin{example}\label{prim2}
Let us consider the 1-dimensional oscillator with respect to the same frame.
The Hamiltonian and Hamiltonian connection of the oscillator read
\ben
&& \cH=\frac12(p^2 + y^2), \nonumber \\
&& \g^y= p, \qquad \g_p = -y. \label{m160}
\een
This is a complete connection. The canonical transformation
\be
y'=y\cos t - p \sin t, \qquad p'=p\cos t + y\sin t 
\ee
brings the connection (\ref{m160}) into zero.
\end{example}

\begin{example}\label{prim3}
Let us consider 1-dimensional motion in a viscous medium 
with respect to the reference 
frame in the previous Examples. It is described by the first order
differential evolution equation
\beq 
y_t=\g^y, \qquad p_t =\g_p, \label{m163}
\eeq
where 
\beq
\g^y =p, \qquad \g_p =-p \label{m162}
\eeq
is a connection on $\Pi$. It is a complete connection, but not
locally Hamiltonian. The coordinate transformation
\beq
y'= y+ p(1-e^t), \qquad p'= pe^t \label{m164}
\eeq
brings the connection (\ref{m162}) into zero so 
that the equations (\ref{m163})
come to the equilibrium equations
\be 
y'_t=0, \qquad p'_t =0.
\ee
However, (\ref{m164}) is not a canonical transformation.
\end{example}

It should be emphasized that, in general, the canonical transformations
introduced above do not preserve the splitting
(\ref{m46}). Consequently, they do not send a Hamiltonian form into 
a Hamiltonian form and do not
maintain the form of the Hamilton equations (\ref{m40}) in general.

At the same time, Proposition \ref{genham} shows that 
canonical morphism send
general Hamiltonian forms to general Hamiltonian forms.

\begin{proposition}\label{genform}
Let $\g$ be a Hamiltonian connection on $\Pi\to X$ and $H$ the corresponding
general Hamiltonian form (see Definition \ref{d5.8}). In virtue of
Proposition \ref{p13.2}, we have $\tau_\g\rfloor\bom=dH$. Let $H'$ be 
another general Hamiltonian form. Then, $\si=H'-H$ is a 1-form on
$\Pi$ such that
\ben
&&d(\si\w dt)= 0, \label{m243} \\
&& \dr_j\si_i -\dr_i\si_j=0, \qquad \dr^j\si^i -\dr^i\si^j=0, \qquad
\dr^j\si_i-\dr_i\si^j=0. \nonumber
\een
\end{proposition}

It follows that, if $\rho$ is a canonical morphism and $H$ is a Hamiltonian
form, then
\be
\rho^*H=H - \si =p_idy^i -(\cH +\si_t)dt -\si_idy^i - \si^idp_i,
\ee
where $\si=\si_tdt+\si_idy^i +\si^idp_i$ is a 1-form on $\Pi$ which 
satisfies (\ref{m243}). Accordingly, the Hamilton equations (\ref{m41a}) --
(\ref{m41b}) are brought into the form
\be
&& y^i_t= \dr^i(\cH+\si_t) +\dr_t\si^i,\\
&& p_{it}= -\dr^i(\cH+\si_t) -\dr_t\si^i.
\ee

\begin{remark} Every general Hamiltonian form is Hamiltonian locally. Every
canonical morphism $\rho$ transforms a Hamiltonian form to a
Hamiltonian form locally since the
condition (\ref{m243}) implies that $\si=fdt + dS$ locally, 
where $f$ and $S$
are local functions on $\Pi$.
\end{remark}

Canonical transformations keep the Hamilton equations if
\beq
\rho^*H=H - dS, \label{m142}
\eeq
where $S$ is a function on $\Pi$ called the {\it generating function}.
In this case, one sais that $\rho^*H$ and $H$ describe the same
mechanical system. The relation (\ref{m142}) can be written as
the Pfaffain equation on the graph $\Pi_\rho\subset\Pi\xx \Pi$ of 
the canonical morphism $\rho$.

In particular, assume that the graph $\Pi_\rho$ is coordinatized by
$(t,y^i,{y'}^i)$. The equality (\ref{m142}) takes the coordinate form
\be
p'_id{y'}^i - p_idy^i + (\cH -\cH')dt = - dS(t,y^i,{y'}^i).
\ee
It leads to the familiar relations
\be
p_i= \frac{\dr S}{\dr y^i}, \quad p'_i= - \frac{\dr S}{\dr {y'}^i}, \quad
\cH'=\cH+ \frac{\dr S}{\dr t}.
\ee

\begin{example} The holonomic coordinate transformations (\ref{m153}) admit
locally the generating function $S(t,{y'}^j,p_i)= y^i(t,{y'}^j)p_i$.
\end{example}

\subsection{Vertical extension of the Hamiltonian formalism}

We now turn to the vertical extension of the time-dependent Hamiltonian
formalism (see Section 4.7) in order to make Hamiltonian forms and 
Hamilton equations invariant under canonical transformations. In case of
symplectic mechanics, the similar extension of the symplectic geometry from
$T^*M$ to $TT^*M$ has been considered in
\cite{abra,tulc74}.

Given a bundle $Y\to X$, let us consider the
Legendre manifold $\Pi_{VY}$ corresponding to the bundle $VY\to X$. It is
isomorphic to the vertical tangent bundle $V\Pi=VV^*Y$ of $\Pi\to X$ (see
Proposition
\ref{vertmom}). We call $\Pi_{VY}$ the {\it vertical phase space} and provide
it with the coordinates
$(x^\la, y^i,  p_i,\dot y^i,\dot p_i)$ of $V\Pi$ (recall the notations
(\ref{m141}) and (\ref{m147})).

The canonical form $\bla$ (\ref{m36}) on $V\Pi$ is the $n=1$ reduction
\beq
\bx{\bla_V=[d\dot p_i\w dy^i +dp_i\w d\dot y^i]\w dt\ot\dr_t}
\label{m144}
\eeq
of the form (\ref{m16}).
The canonical form $\bom$ (\ref{m66}) on $V\Pi$ reads
\beq
\bx{\bom_V=[d\dot p_i\w dy^i +dp_i\w d\dot y^i]\w dt.}
\label{m145}
\eeq

With the canonical form (\ref{m145}), the vertical  phase space
$V\Pi$ can be equipped with the canonical Poisson structure (\ref{m72})
given by the bracket
\beq
\bx{\{f,g\}_{VV} =\dot\dr^if\dr_ig +\dr^if\dot\dr_ig -\dr^ig\dot\dr_if
-\dot\dr^ig\dr_if.} \label{m146}
\eeq

The notions of Hamiltonian connection, 
Hamiltonian vector field, horizontal Hamiltonian vector field and 
Hamiltonian form on $V\Pi$ are the straightforward generalization of those
in Section 5.1.

Every Hamiltonian form $H_V$ on $V\Pi$ admits the splitting
\beq
\bx{H_V=\dot p_idy^i -\dot y^idp_i -\cH_V, \quad \cH_V =(\dot p_i\wt\g^i 
-\dot y^i\wt\g_i +\wt\cH_V)dt,} \label{m148}
\eeq
where $\g$ is a connection on $\Pi\to X$. The corresponding Hamilton
equations (\ref{m40}) takes the form
\bea
&& \g^i=\dot\dr^i\cH_V, \label{m149a}\\
&& \g_i=-\dot\dr_i\cH_V,\label{m149b}\\
&& \dot \g^i=\dr^i\cH_V,\label{m149c}\\
&& \dot \g_i=-\dr_i\cH_V.\label{m149d}
\eea

In particular, the vertical lift $V\wt\g$ (\ref{43}) of a Hamiltonian
connection $\wt\g$ associated with a Hamiltonian $\cH$ on $\Pi$ is the
Hamiltonian connection associated with the Hamiltonian 
\be
\cH_V=\dr_V\cH=(\dot y^i\dr_i +\dot p_i\dr^i)\cH
\ee
on $V\Pi$ (see Proposition \ref{p01}). In this case, the Hamilton equations
(\ref{m149a}) and (\ref{m149b}) are 
exactly the Hamilton equations (\ref{m40})
for the Hamiltonian connection $\wt\g$. 

Let us consider the canonical coordinate 
transformations of the Legendre bundle
$\Pi\to X$ and the induced (holonomic) coordinate transformations 
\beq
\dot p'_i=\dr_Vp'_i, \qquad \dot{y'}^i=\dr_V{y'}^i \label{m151}
\eeq
of $V\Pi$.
It is readily observed that they are also canonical transformations for
the canonical forms (\ref{m144}), (\ref{m145}). They are linear in the
coordinates $\dot p_i$, $\dot y^i$ and, obviously, do not exhaust all
canonical transformations of $V\Pi$. 
These transformations maintain the Poisson bracket (\ref{m146}).
The splitting (\ref{m148}) of a 
Hamiltonian $\cH_V$ and the Hamilton equations
(\ref{m149a}) -- (\ref{m149d}) also are invariant under the canonical
transformations (\ref{m151}). We have
\be
H_V=\dot p_idy^i -\dot y^idp_i -\cH_V=\dot p'_id{y'}^i -\dot{y'}^idp'_i
-\cH'_V,
\ee
where
\beq
\cH'_V=\cH_V - (\dr_Vp'_i\dr_t{y'}^i - \dr_V{y'}^i\dr_tp'_i). \label{m150}
\eeq
At the same time, if
$\cH_V=\dr_V\cH$, where $\cH$ is a Hamiltonian on $\Pi$, the Hamiltonian
$\cH'_V$ (\ref{m150}) fails to represent the derivative $\dr_V$ of some
Hamiltonian on $\Pi$ in general.

\begin{proposition}\label{can4}
Every connection $\wt\g$ on the Legendre bundle $\Pi$ gives rise to the
Hamiltonian connection on $V\Pi$. 
\end{proposition}

\begin{proof}
Let us consider the Hamiltonian form
\be
H_V=\dot p_i(dy^i-\wt\g^idt) -\dot y^i(dp_i -\wt\g_idt)  = 
\dot p_idy^i -\dot y^idp_i -(\dot p_i\wt\g^i -\dot y^i\wt\g_i)dt
\ee
The corresponding Hamiltonian connection on $V\Pi$ is given by the
Hamilton equations (\ref{m149a}) -- (\ref{m149d}) which take the form
\beq
\g^i=\wt\g^i, \quad \g_i=\wt\g_i, \quad \dot\g^i =\dot p_j\dr^i\wt\g^j-\dot
y^j\dr^i\wt\g_j, \quad \dot\g_i = -\dot p_j\dr_i\wt\g^j + \dot
y^j\dr_i\wt\g_j. \label{m165}
\eeq
In particular, if $\wt\g$ is a Hamiltonian connection on $\Pi$, the
Hamiltonian connection (\ref{m165}) coincides with the vertical
 connection $V\g$ (\ref{43}).
\end{proof}

It follows that every first order evolution equations (\ref{m170}) 
on the Legendra bundle $\Pi$ can be written as a part 
(\ref{m149a}), (\ref{m149b}) of the Hamilton equations on $V\Pi$.

\begin{example}\label{prim4}
The 1-dimensional motion in a viscous medium in Example \ref{prim3} is
described by the Hamiltonian  $\cH_V= p(\dot p  +\dot y)$ on $V\Pi$.
\end{example}

\subsection{Reference frames}

The form of the Hamiltonian evolution equation 
(\ref{m96}) is maintained under
canonical transformations of $\Pi$ when
\be
\wt\cH'(t,{y'}^j,p'_j)= \wt\cH(t,y^j,p_j), \qquad
(\dr_t +\g_0{}^i\dr_i +\g_{0i}\dr^i)f =
(\dr_t +\g'_0{}^i\dr'_i +\g'_{0i}{\dr'}^i)f'.
\ee

In virtue of Corollary \ref{c13.5}, we can make locally the Hamiltonian
connection $\g_0$ equal to zero by  canonical coordinate transformations and
can bring the Hamiltonian evolution equation (\ref{m96}) into the familiar
Poisson bracket form
\beq
d_{Ht}f= \dr_tf +\{\wt\cH,f\}_V. \label{m152}
\eeq
In virtue of Proposition \ref{p13.5}, we can get this form of the
Hamiltonian evolution equation with respect to the global trivialization of
$\Pi$ if $X=\R$ and the Hamiltonian connection
$\g_0$ in the splitting (\ref{m95}) is complete.

In particular, let $\G$ be a complete connection on the 
bundle $Y\to\R$ associated with some trivialization (\ref{m33}) of $Y$.
Then, the connection
$\wt\G=V^*\G$ (\ref{m38}) is a complete Hamiltonian connection on $\Pi$ 
associated with the corresponding trivialization (\ref{m35}) of $\Pi$.
 It follows that we can utilize the
covector lift $\g_0=\wt\G$ of a complete connection $\G$ on $Y$ in order to
bring the Hamiltonian evolution equation (\ref{m96}) into the Poisson bracket
form (\ref{m152}). 

\begin{definition}\label{refr2}
We say that a complete connection $\G$ on $\Y$ describes a {\it reference
frame} in time-dependent mechanics.
\end{definition}

Indeed, the difference of $\G'-\G=u dt$ defines a vertical vector
field $u$ on $Y$ which characterizes the relative 
velocities between reference
frames $\G'$ and $\G$. Accordingly, one can think of
\be
\dot y^i\circ D_\G =y^i_t-\G^i
\ee
as being the relative velocities of a mechanical system with respect to the
reference frame $\G$.
By Definition \ref{refr2}, there is the 1:1 correspondence between
reference frames and trivializations of $Y\to \R$. 

\begin{imp}One can say that a reference
frame provides a splitting between the time and the other coordinates of a
mechanical system.\end{imp}

\begin{remark}\label{refr3}
Every connection $\G$ on $\Y$
defines a local reference frame. In virtue of Proposition \ref{p5.12}, every
Hamiltonian form $H$ on
$\Pi$ defines the connection $\G_H$ on $\Y$ and,
consequently, the corresponding local reference frame. We call $\G_H$ the
(local) {\it proper reference frame}. 
With respect to this reference frame, the
Hamilton equations (\ref{m41a}) takes the form
\beq
y^i_t-\G_H=\dr^i\wt\cH. \label{m154}
\eeq
One can think of these equations as 
being the relations between the canonical
momenta $p_i$ and the velocities $y^i_t-\G^i_H$ relative
to the proper reference frame. In accordance with the definition of $\G_H$,
this relation implies that the null momenta corresponds to the null 
velocities (\ref{m154}).
\end{remark}

\newpage

\section{Lagrangian mechanics}

We aim to investigate the relations between Lagrangian and 
Hamiltonian formulations of time-dependent mechanics. From the mathematical
point of view, these formulations are not equivalent 
in case of degenerate Lagrangians. From physical 
viewpoint, velocities are physical observables in classical
mechanics, whereas momenta are physical observables in quantum mechanics.

Given a bundle $\Y$ over a 1-dimensional base $X$, 
the Lagrangian mechanics of
sections of
$\Y$ is formulated on the configuration space $J^1Y$ coordinatized by
$(t,y^i, y^i_t)$. A Lagrangian on $J^1Y$ reads $L=\cL dt$. Also recall the
notation
\be
\pi_i=\dr_i^t\cL, \qquad \pi_{ij}=\dr^t_j\dr_i^t\cL, \qquad \wt\pi=\cL-\pi_i
y^i_t.
\ee

\subsection{Poisson structure}

In contrast with Hamiltonian mechanics, the configuration space $J^1Y$ of
Lagrangian mechanics possesses no canonical Poisson structure. 

Given a Lagrangian $L$ on $J^1Y$, 
the pullbacks  of the canonical forms on
$V^*Y$ and $T^*Y$ are defined by the Legendre morphism
$\wh L: J^1Y\to V^*Y$ (\ref{m11}) and the morphism $\wh\Xi_L: J^1Y\to T^*Y$
(\ref{N42}) on $J^1Y$.

Let $\bom$ be the canonical 3-form (\ref{m66}) on $V^*Y$. Its pullback by the
Legendre morphism $\wh L$ reads
\be
\bom_L=\wh L^*\bom =d\pi_i\w dy^i\w dt= (\pi_{ij}dy^j_t\w dy^i + \dr_j\pi_i
dy^j\w dy^i)\w dt. 
\ee
Using $\bom_L$, every vertical vector field 
$\vt=\vt^i\dr_i +\dot\vt^i\dr_i^t$
on
$J^1Y\to X$ is sent to the 2-form
\be
\vt\rfloor\bom_L= \{[\dot\vt^j\pi_{ji} +\vt^j(\dr_j\pi_i -\dr_i\pi_j)]dy^i-
\vt^i\pi_{ji}dy^j_t\}\w dt.
\ee 
If the Lagrangian $L$ is regular
($\det\pi_{ij}\neq 0$), the above map is a bijection. 
Indeed, given any 2-form
$\f=(\f_idy^i +\dot\f_idy^i_t)\w dt$, the algebraic equations 
\be
 \dot\vt^j\pi_{ji} +\vt^j(\dr_j\pi_i -\dr_i\pi_j)=\f_i, \qquad
 -\vt^i\pi_{ji}=\dot\f_j
\ee
have the unique solution
\be
\vt^i=-(\pi^{-1})^{ij}\dot\f_j, \qquad \dot\vt^j=(\pi^{-1})^{ji} [\f_i +
(\pi^{-1})^{kn}\dot\f_n(\dr_k\pi_i -\dr_i\pi_k)].
\ee 
In particular, every function $f$ on $J^1Y$ defines a vertical vector field
\beq
\vt_f=-(\pi^{-1})^{ij}\dr_j^tf\dr_i + (\pi^{-1})^{ji} [\dr_if +
(\pi^{-1})^{kn}\dr_n^tf(\dr_k\pi_i -\dr_i\pi_k)]\dr_j^t. \label{m125}
\eeq
Following the relation (\ref{m124}),  
 one can introduce the Poisson structure on the space $S(J^1Y)$ of
functions on $J^1Y$. It is given by the bracket
\ben
&&\vt_g\rfloor\vt_f\rfloor\bom_L=\{f,g\}_Ldt,  \nonumber \\
&&\{f,g\}_L= [(\pi^{-1})^{ij} +(\dr_n\pi_k
-\dr_k\pi_n)(\pi^{-1})^{ki}(\pi^{-1})^{nj}](\dr_i^tf\dr_jg-\dr_i^tg\dr_jf) +
\label{m121} \\
&&\qquad (\dr_n\pi_k
-\dr_k\pi_n)(\pi^{-1})^{ki}(\pi^{-1})^{nj}\dr_i^tf\dr_j^tg. \nonumber
\een
The vertical vector field 
$\vt_f$ (\ref{m125}) is the Hamiltonian vector field
of the function $f$ with respect to this Poisson structure.

In particular, if the Lagrangian $L$ is hyperregular, that is, the
Legendre morphism $\wh L$ is diffeomorphism, the Poisson structure
(\ref{m121}) is obviously isomorphic to the Poisson structure (\ref{m72}) on
the  phase space $\Pi=V^*Y$. 

The Poisson structure (\ref{m121}) defines the corresponding symplectic
foliation on $J^1Y$ which consists with the fibration $J^1Y\to X$. The
symplectic form on the leaf $J^1_xY$ of this foliation is 
$\Om_x= d\pi_i\w dy^i$ \cite{vais95}.

\subsection{Spray-like equations}

In the framework of Hamiltonian mechanics above, 
we have shown that the choice
of a trivialization
$Y\simeq X\xx M$ corresponds to the choice of a certain reference frame. We
here illustrate this fact in case of evolution equations on the configuration
space.
We consider second order evolution equations which are not necessarily
of Lagrangian type. 

Let us recall the notion of spray in autonomous mechanics.
Let $M$ be a manifold coordinatized by $(\ry^i)$ and 
\be
K=d\ry^i\ot(\dr_i - K^k{}_{ji}\dot\ry^j\dot\dr_k) 
\ee
a linear connection (\ref{408}) on $TM$. It yields the vector field 
\beq
\dot\ry^i\dr_i\rfloor K(\ry,\dot\ry)= \dot\ry^i(\dr_i
- K^k{}_{ji}\dot\ry^j\dot\dr_k) \label{m172}
\eeq
on $TM$ which is called the {\it geodesic spray}. The equations
\be
\frac{d\ry^i}{dt}= \dot\ry^i, \qquad
\frac{d\dot\ry^i}{dt} =- K^i{}_{ji}\dot\ry^j\dot\ry^i 
\ee
for integral curves of the spray (\ref{m172}) are second order
differential equations whose solutions are geodesics of the connection $K$. 

We aim to discover similar spray-like equations in time-dependent
mechanics.

Given a bundle $Y\to X$ over a 1-dimensional base, let 
\beq
A=dt\ot(\dr_t + A^i\dr_i^t) + dy^j\ot(\dr_j + A^i_j\dr_i^t) \label{m173}
\eeq
be a connection on the jet bundle $J^1Y\to Y$. It has the transformation law
\ben
&& {A'}^i_k = (\frac{\dr{y'}^i}{\dr y^j}A^j_n
+\frac{\dr{y'}^i_t}{\dr y^n})\frac{\dr y^n}{\dr{y'}^k}, \label{m175}\\
&&{A'}^i =(\frac{\dr{y'}^i}{\dr y^j}A^j + \frac{\dr{y'}^i_t}{\dr t})
+ (\frac{\dr{y'}^i}{\dr y^j}A^j_k +\frac{\dr{y'}^i_t}{\dr y^k})
\frac{\dr y^k}{\dr t} = \frac{\dr{y'}^i_t}{\dr y^j_t}A^j +
\frac{\dr{y'}^i_t}{\dr t} - {A'}^i_k\frac{\dr{y'}^k}{\dr t}. \nonumber
\een

\begin{definition}\label{spray}
Given a connection $A$ (\ref{m173}) on $J^1Y\to Y$, by a {\it second order
evolution equation} on $Y$ is meant the restriction of the kernel $\Ker\wt
D\subset J^1J^1Y$ of the vertical covariant differential
$\wt D$ (\ref{7.10}) to $J^2Y$. This is given by the coordinate relation
\beq
y^i_{tt}=  A^i +A^i_{tj}y^j. \label{m174}
\eeq
\end{definition}

A glance at the expression (\ref{m174}) shows that different connections
(\ref{m173}) can lead to the same evolution equation.

\begin{remark}
Every connection (\ref{m173}) on $J^1Y\to Y$ generates the connection
\beq
\g=dt\ot(\dr_t + y^i_t\dr_i + (A^i +A^i_j y^j_t)\dr_i^t) \label{m213}
\eeq
on $J^1Y\to X$. The horizontal lift of the vector field $\dr_t$ on $X$ onto
$J^1Y$ by means of the connection (\ref{m213}) reads
\be
\dr_t + y^i_t\dr_i + ( A^i +A^i_j y^j_t)\dr_i^t. 
\ee
The integral curves of this vector field are the generalized solutions of
the evolution equations (\ref{m174}). Conversely, second order evolution
equation can is often defined as the equation 
\be
d_t y^i=y^i_t, \qquad d_ty^i_t=\xi^i,
\ee
for an integral curve of a vector field
\be
\dr_t + y^i_t\dr_i + \xi^i\dr_i^t 
\ee
on $J^1Y$. Every such a vector field defines a connection
\be
A^i_j=\frac12\frac{\dr\xi^I}{\dr y^j_t}, \qquad A^i= \xi^i -A^i_j
\ee
on $J^1Y\to Y$ \cite{giac93} which lead to the same evolution equation in
accordance with Definition \ref{spray}.
\end{remark}

In particular, let $\Y$ be a trivializable bundle and $Y\simeq X\times M$ its
trivialization with the coordinates $(t,\ry^i)$ 
whose transition functions are
independent on $t$. We have the corresponding trivialization
$J^1Y\simeq X\times TM$ with the coordinates $(t,\ry^i,\dot\ry^i)$, where
$\dot\ry^i$ are  holonomic coordinates of $TM$.  With respect to
these coordinates, the transformation law (\ref{m175}) of the connection
(\ref{m173}) reads
\beq
{A'}^i =\frac{\dr\ry'^i}{\dr\ry^j}A^j \qquad
{A'}^i_k = (\frac{\dr\ry'^i}{\dr \ry^j}A^j_n
+\frac{\dr \dot \ry'^i}{\dr \ry^n})\frac{\dr\ry^n}{\dr{\ry'}^k}.
\label{m176}
\eeq

A glance at the expression (\ref{m176}) 
shows that, given a trivialization of
$\Y$, a connection on $J^1Y\to Y$ defines a vertical time-dependent vector
field
$A^i_t$ on $TM$ and a time-dependent connection on $TM\to M$.
The converse procedure enables us to discover a spray-like equation on $Y$.

Let $K^i_k(\ry^j, \dot\ry^j)$ be a connection (e.g., a linear
connection) on $TM\to M$.  
Given the above-mentioned trivialization of $Y$, the connection
$K$ defines the connection $A$ on
$J^1Y\to Y$ by the coordinate relations
\beq
A^i=0, \qquad A^i_k= K^i_k. \label{m177}
\eeq
Owing to the transformation law (\ref{m175}), we can write the
connection (\ref{m177}) with respect to arbitrary bundle coordinates
$(t,y^i)$. It reads
\ben
&&A^i_k=[\frac{\dr y^i}{\dr\ry^j}K^j_n(\ry^j(y^i), \dot\ry^j(y^i,y^i_t))
+\frac{\dr^2 y^i}{\dr\ry^n\dr\ry^j} \dot\ry^j +
\frac{\dr\G^i}{\dr\ry^n}]\dr_k\ry^n,  \label{m178}\\
&&A^i = \dr_t\G^i +\dr_j\G^iy^j_t - A^i_k\G^k, \nonumber
\een
where $\G^i=\dr_ty^i$ is the connection on $\Y$ which
corresponds to the initial trivialization of $Y$, that is, $\G=0$
relative to the coordinates $(t,\ry^i)$.

\begin{remark}
Given the connection $A$ (\ref{m177}) on $J^1Y\to Y$ and the
above-mentioned connection 
$\G$ on $Y\to X$, the corresponding 
composite connection (\ref{1.39}) consists
with the jet lift $J\G$ (\ref{59}) of $\G$ onto $J^1Y\to Y$. We have
\beq
J\G=dt\ot(\dr_t + \G^i\dr_i +d_t\G^i\dr_i^t). \label{m190}
\eeq
\end{remark}

The evolution equation (\ref{m174}) with respect to the connection
(\ref{m178}) reads
\beq
y^i_{tt}=\dr_t\G^i +y^j_t\dr_j\G^i + A^i_k(y^k_t-\G^k). \label{m179}
\eeq

Given a reference frame $Y\simeq X\xx M$ with coordinates $(t,\ry^i)$, let
$K^i_j=0$. In accordance with (\ref{m177}), this choice leads to the free
motion equation 
\beq
\ddot\ry^i=0. \label{m180}
\eeq
With respect to arbitrary bundle coordinates $(t,y^i)$, this equation reads
\beq
y^i_{tt}=\dr_t\G^i -\G^j\dr_j\G^i+y^k_t(2\dr_k\G^i+
\frac{\dr y^i}{\dr\ry^j}\frac{\dr\ry^j}{\dr y^m\dr y^k}\G^m)
-\frac{\dr y^i}{\dr\ry^j}\frac{\dr\ry^j}{\dr y^m\dr y^k}y^k_ty^m_t.
\label{m188}
\eeq
One can treat the right side of this equation as the general expression for
inertial forces. 

\begin{imp} Such kind of terms in spray-like 
evolution equations can be always
excluded by the choice of a reference frame. \end{imp}

\begin{remark} The equations (\ref{m180}) are obviously of Lagrangian
type, but not the equations (\ref{m188}). At the same time, the equations
(\ref{m188}) are equivalent to the 
Euler--Lagrange equations of the Lagrangian
\be
\cL=\frac12\dl_{ab}\frac{\dr y^a}{\dr{y'}^i}
\frac{\dr y^b}{\dr{y'}^j}({y'}^i_t
-\G^i) ({y'}^j_t -\G^j).
\ee
\end{remark}

\begin{example}\label{rot}
Let us consider a free point on a plain. 
Let the splitting  $Y=\R\times\R^2$
with coordinates $(t,\ry^1,\ry^2)$ corresponds to an inertial reference
frame. Let the connection $K$ on the bundle $T\R^2$ be equal to zero.
Let us consider the rotating reference frame with the coordinates
\be
y^1=\ry^1\cos wt -\ry^2\sin wt, \qquad y^2=\ry^2\cos wt +\ry^1\sin wt.
\ee
With respect to this reference frame, the equation (\ref{m188}) reads
\beq
y^i_{tt}=\dr_t\G^i +2y^j_t\dr_j\G^i -\G^j\dr_j\G^i, \label{m182}
\eeq
where 
\beq
\G^1=\dr_t y^1=-wy^2, \qquad \G^2=\dr_t y^2=wy^1. \label{m183}
\eeq
Substituting (\ref{m183}) into (\ref{m182}), we find
\be
y^1_{tt} =w^2y^1- 2wy^2, \qquad y^2_{tt} =w^2y^2 + 2wy^1.
\ee
\end{example}

This Example shows that, on physical level, we can treat $y^i_t$ in the
evolution equation (\ref{m179}) as the velocities relative to the
(local) reference frame given by the connection on $\Y$ which vanishes with
respect to these coordinates (see Section 5.6).

At the same time, the evolution equation (\ref{m179}) is brought into the
spray-like form
\be
&& d_t\dot y^i= {K'}^i_k\dot y^k +\dr_k\G^i\dot y^k, \\
&&  {K'}^i_k= 
\frac{\dr y^i}{\dr\ry^j}\frac{\dr\ry^n}{\dr y^k}
K^j_n(\ry^j(y^i), \dot\ry^j(y^i,\dot y^i)) -
\frac{\dr y^i}{\dr\ry^j}\frac{\dr\ry^j}{\dr y^m\dr y^k}\dot y^m, \\ 
&& \dot y^i =y^i_t -\G^i, \qquad \dot\ry^i=\frac{\dr\ry^i}{\dr y^k}\dot y^k,
\ee
where $\dot y^i\dr_i$ can be treated the relative velocities with respect to
the initial reference frame $\G$ which are written with respect to the
coordinates $y^i$. 

In particular, if $K^i_k=- K^i{}_{mk}\dot\ry^m$ is a
linear connection on $TM$, we have
\be
&& d_t\dot y^i= -{K'}^i{}_{mk}\dot y^m\dot y^k +\dr_k\G^i\dot y^k, \\
&&  {K'}^i_{mk}= 
\frac{\dr y^i}{\dr\ry^j}\frac{\dr\ry^n}{\dr y^k}\frac{\dr\ry^l}{\dr y^m}
K^j{}_{lk} -
\frac{\dr y^i}{\dr\ry^j}\frac{\dr\ry^j}{\dr y^m\dr y^k}\dot y^m.
\ee

\subsection{Hamiltonian and Lagrangian formalisms}

According to Section 4.6, we establish the relations
between Lagrangian and Hamiltonian formalisms for time-dependent mechanical
systems.

Let $\Y$ be a bundle over a 1-dimensional base, $\Pi=V^*Y$ the phase space 
and $J^1Y$ the configuration space.

\begin{definition}
A Hamiltonian 
$\cH$ on $\Pi$ is said to be associated with a Lagrangian $L$ on $J^1Y$ if it
obeys the conditions
\bea 
&&\bx{\wh L\circ\wh H|_Q=\Id_Q, \quad 
p_i=\dr^t_i \cL(t, y^j, 
\dr^j\cH(p)), \quad p\in Q=\wh L(J^1Y)} \label{m191a}\\
&& \bx{ H_{\wh H}-H=L\circ\wh H, \qquad \cL(t, y^j, 
\dr^j\cH)\equiv p_i\dr^i\cH -\cH.} \label{m191b}
\eea
\end{definition}

Also the relation
\beq
\bx{\dr_i\cH +\dr_i\cL=0} \label{m192}
\eeq
takes place on $Q$.

If a Lagrangian $L$ is hyperregular, there exists an unique
Hamiltonian associated with $L$.
If a Lagrangian $L$ is degenerate, 
different Hamiltonians or no Hamiltonian at all may be associated with
$L$ in general.

\begin{proposition}\label{almreg} Let a Lagrangian $L$ be {\it almost
regular}, that is, the Lagrangian constraint space
$Q$ is an imbedded submanifold of $\Pi$ and the
Legendre morphism $\wh L:J^1Y\to Q$ is a submersion. Then, each point of $Q$
has an open neighborhood on which there exists a local
Hamiltonian form associated with $L$ \cite{sard95,zakh}.
\end{proposition}
 
\begin{example}
Let $Y$ be the bundle $\R^2\to\R$ coordinatized
by $(t,y)$. Consider the Lagrangian $\cL=\exp y_t$.
It is regular and semiregular, but not hyperregular. The corresponding
Legendre morphism reads $p\circ\wh L=\exp y_t$.
The image $Q$ of the configuration space under this morphism
is given by the coordinate relation $p>0$. It is an open
subbundle of the Legendre bundle.
On $Q$, we have the associated Hamiltonian 
\be
\cH=p(\ln p -1)
\ee
which however can not be extended to $\Pi$.
\end{example}

All Hamiltonian forms
associated with a semiregular Lagrangian $L$ coincide with each
other on the Lagrangian constraint space $Q$, and the Poincar\'e--Cartan form
$\Xi_L$ (\ref{303}) is the pullback 
\be
\bx{\pi_iy^i_t-\cL\equiv \cH(t,y^i,\pi_i),} 
\ee
of such a Hamiltonian form $H$ by the Legendre morphism $\wh L$.

\begin{example}
Let $Y$ be the bundle $\R^3\to\R$ coordinatized
by $(t,y^1,y^2)$. Consider the Lagrangian
\beq
\cL=\frac12(y^1_t)^2. \label{200}
\eeq
It is semiregular. The associated Legendre morphism reads
\be
p_1\circ\wh L=y^1_t, \qquad  p_2\circ\wh L= 0. 
\ee
The corresponding constraint space $Q$ consists of points with
the coordinate $p_2=0$. The Hamiltonians associated with
the Lagrangian (\ref{200}) are given by the expression
\beq
\cH=\frac12 (p_1)^2 + c(t,y)p_2, \label{m196}
\eeq
where $c$ is arbitrary function on $Y$. They coincide with each other on $Q$.
\end{example}

\begin{corollary}\label{j3}
In accordance with Proposition \ref{pHL},
if $\cH$ is a Hamiltonian associated with a semiregular
Lagrangian $L$, every solution of the corresponding Hamilton equations
which  lives on the Lagrangian constraint space $Q$
yields a solution of the Euler--Lagrange equations for $L$.
At the same time, to exaust all solutions of the Euler--Lagrange
equations, one must consider a complete family (if it exists) of
Hamiltonians associated with $L$.
\end{corollary}

For instance, the Hamiltonians (\ref{m196}) associated with the
Lagrangian (\ref{200}) constitute a complete family. 

\begin{proposition}\label{comfam} If a 
Lagrangian $L$ is semiregular and almost
regular, then every point of $Q$ has a neighborhood on which there exists a
complete family of local Hamiltonians associated with $L$ \cite{sard95,zakh}.
\end{proposition} 

The following example shows
that a complete family of associated Hamiltonians
may exist even if a Lagrangian is neither semiregular nor almost regular.

\begin{example}
Let $Y$ be the bundle $\R^2\to\R$. Let us consider the Lagrangian 
\be
\cL=\frac13(y_t)^3.
\ee
The associated Legendre morphism reads
\beq
p\circ\wh L=(y_t)^2. \label{204}
\eeq
The corresponding constraint space $Q$ is given by the coordinate
relation $p\geq 0$. It is not even a submanifold of $\Pi$.
There exist two associated Hamiltonians
\be
\cH_+=\frac23 p^{\frac32},\qquad \cH_-=-\frac23 p^{\frac32} 
\ee
which are defined only on the constraint space $Q$. They correspond to
different solutions
\be
y_t=\sqrt{p}, \qquad  y_t=-\sqrt{p}
\ee
of the equation (\ref{204}) and constitute a complete family.
\end{example}

\subsection{Quadratic Lagrangians and Hamiltonians}

As an important illustration of Proposition \ref{comfam}, 
let us describe the
complete families of Hamiltonians associated with almost regular quadratic
Lagrangians.

\begin{remark}
Since Hamiltonians in time-dependent mechanics are not functions on a phase
space, we can not apply to them the well-known analysis of the normal forms
\cite{brun} (e.g. quadratic Hamiltonians in sypmlectic mechanics
\cite{arno}).
\end{remark}

Let us consider a quadratic Lagrangian 
\beq
\cL=\frac12 a_{ij}(y)y^i_t y^j_t +
b_i(y)y^i_t + c(y) \label{N12}
\eeq
on $J^1Y$, where $a$, $b$ and $c$ are local functions on $Y$ with the
corresponding transformation laws. The associated Legendre morphism reads
\beq
p_i\circ\wh L= a_{ij} y^j_t +b_i. \label{N13}
\eeq
It is easily observed that the Lagrangian (\ref{N12}) is semiregular. 

The Legendre morphism (\ref{N13})
is an affine morphism over $Y$. The 
corresponding linear morphism over $Y$ is
\be
\ol L: VY\to V^*Y, \qquad p_i\circ\ol L=a_{ij}\dot y^j,
\ee
where $\dot y^j$ are bundle coordinates of the vector
bundle (\ref{23}). In particular, if $L$ is regular, the morphism $\ol L$
defines a nondegenerate fibre metric on $VY$.

Let us assume 
that the Lagrangian is almost regular (see Proposition \ref{almreg}) and that
the Lagrangian constraint space $Q$
defined by the Legendre morphism (\ref{N13})
contains the image of the 
zero section $\wh 0(Y)$ of the Legendre bundle $\Pi\to Y$.
It is immediately observed that
$\Ker\wh L= \wh L^{-1}(\wh 0(Y))$ is an affine subbundle of the jet bundle
$J^1Y\to Y$.

The following two ingredients in our constrauction play a prominent role.

(i) There exists a connection $\G$ on $Y\to X$
which takes its values into $\Ker\wh L$:
\beq
\G: Y\to \Ker\wh L, \qquad \bx{a_{ij}\G^j + b_i =0.} \label{250}
\eeq
With this connection, the Lagrangian 
(\ref{N12}) can be brought into the form
\be
\cL=\frac12 a_{ij}(y^i_t -\G^i)(y^j_t -\G^j) + c'. 
\ee
For instance, if it is regular, the connection (\ref{250}) is unique.

(ii) There exists a linear morphism
\beq
\si: V^*Y\to VY, \qquad \dot y^i\circ\si =\si^{ij}p_j \label{N17}
\eeq
such that
\be
\ol L\circ\si\mid_Q= \Id_Q, \qquad
\bx{a_{ij}\si^{jk} a_{kb}=a_{ib}.} 
\ee
Then, the jet bundle $J^1Y\to Y$ has the splitting
\beq
\bx{J^1Y=\Ker\wh L\op\oplus_Y{\rm Im}\si,} \qquad y^i_t= [y^i_t
-\si^{ik}(a_{kj}y^j_t + b_k)]+
[\si^{ik}(a_{kj}y^j_t + b_k)]. \label{N18}
\eeq
If the Lagrangian (\ref{N12}) is regular, the morphism (\ref{N17})
is determined uniquely.

Given the morphism $\si$ (\ref{N17}) and the connection $\G$
(\ref{250}), let us consider  the Hamiltonian form 
\beq
H= p_idy^i - [\G^i
(p_i-\frac12 b_i) +\frac12 \si^{ij}p_ip_j-c]dt. \label{N22}
\eeq

\begin{proposition} The Hamiltonian form (\ref{N22})
is associated with the Lagrangian (\ref{N12}). The family of these forms
parameterized by the connections (\ref{250}) constitute a complete
family.
\end{proposition}

Given the Hamiltonian (\ref{N22}),
let us consider the Hamilton equations (\ref{m41a}) for
sections $r$ of the bundle $\Pi\to X$. They read
\ben
&&J^1s= (\G+\si)\circ r, \qquad s=\pi_{\Pi Y}\circ r, \label{N29}\\
&&d_t r^i=\G^i + \si^{ij}r_j. \nonumber
\een

With the splitting (\ref{N18}), we have the following surjections
\be
&& {\cal S}:=\pr_1: J^1Y\to\Ker\wh L, \qquad {\cal S}: y^i_t\to y^i_t
-\si^{ik} (a_{kj}y^j_t + b_k), \\
&&\cF:=\pr_2: J^1Y\to{\rm Im}\si, \qquad
\cF=\si\circ\wh L:
y^i_t\to \si^{ik} (a_{kj}y^j_t + b_k).
\ee
With respect to these surjections, the Hamilton equations (\ref{N29})
break into two parts
\ben
&&{\cal S}\circ J^1s=\G\circ s, \qquad
d_t r^i -
\si^{ik} (a_{kj}d_t r^j + b_k)= \G^i, \label{N23}\\
&& \cF \circ J^1s=\si\circ r, \qquad
\si^{ik} (a_{kj}d_t r^j + b_k)= \si^{ik}r_k. \nonumber
\een
The Hamilton equations (\ref{N23}) are independent of canonical momenta
and play the role of constraints.

It should be noted that the Hamiltonian (\ref{N22}) differ from each
other only in connections $\G$ (\ref{250}) which lead to the different
constraints (\ref{N23}).

\begin{remark} We observe that a mechanical system described by a degenerate
Lagrangian $L$ appears a multi-Hamiltonian 
constrained system in the framework
of the Hamiltonian formalism. In the spirit of the
well-known Gotay algoritm in autonomous mechanics
\cite{berg,got78}, the Lagrangian constraint space $Q$ can be called the
primary constraint space. To properly apply 
this algoritm, however, one has to
consider each Hamiltonian of a complete family of Hamiltonians associated
with
$L$. If $L$ is semiregular, all these Hamiltonians coincide with each other
on $Q$, but not the horizontal Hamiltonian vector fields (\ref{m57}). A
different way is to utilize the Gotay algoritm in the framework of the
Lagrangian formalism \cite{leon93,rana}. One can investigate  also the
conditions of formal integrability \cite{gols,kras,pomm} of the Hamilton
equations. Given a Hamiltonian associated with $L$, the corresponding
Hamilton equations fail to satisfy these conditions at all poits of $Q$.
\end{remark}

The relations between Lagrangian and Hamiltonian formalisms described 
above are broken under canonical transformations if the
transition functions $y^i\to {y'}^i$ depend on momenta. In the next Section,
we overcome this difficulty.

\subsection{The unified Lagrangian and Hamiltonian formalism}

In case of a 1-dimensional base $X$, we can generalize the construction given
in Remark \ref{joint} as follows.

Given a bundle $Y\to X$, let $V^*J^1Y$ be the
vertical cotangent bundle of $J^1Y\to X$ coordinatized by 
$(t,y^i,y^i_t, \dot
y_i, \dot y_i^t)$ and
$J^1V^*Y$ the jet manifold of $V^*Y\to X$ 
coordinatized by $(t,y^i,p_i,y^i_t,
p_{it})$. 

\begin{lemma}
There is the isomorphism
\beq
\bpi=V^*J^1Y= J^1V^*Y, \qquad \dot y_i \llra p_{it}, \qquad \dot y_i^t\llra
p_i, \label{m197}
\eeq
over $J^1Y$.
\end{lemma}

\begin{proof}
The isomorphism (\ref{m197}) is proved by compairing the 
transition functions of the coordinates $(\dot y_i, \dot y_i^t)$ and $(p_i,
p_{it})$. 
\end{proof}

Due to the isomorphism (\ref{m197}), one can think of 
$\bpi$ as being both the Legendre bundle (phase space) over the configuration
space $J^1Y$ and the configuration space over the phase space
$\Pi$. Hence, the space $\bpi$ can be utilized as the unified
configuration and phase space of the joint Lagrangian and Hamiltonian
formalism. This space is coordinatized by $(t,y^i,y^i_t, p_{it},p_i)$,
 where $(y^i,p_{it})$ and $(y^i_t,p_i)$ are 
canonically conjugate pairs. The space $\bpi$ is equipped with the canonical
form (\ref{m36}) given by the coordinate form
\be
\bx{\bla=(dp_{it}\w dy^i +dp_i\w dy^i_t)\w dt\ot\dr_t} 
\ee
and with the canonical form (\ref{m66}) which reads
\beq
\bx{\bom=(dp_{it}\w dy^i +dp_i\w dy^i_t)\w dt=d_t(dp_i\w dy^i\w dt).}
\label{m200}
\eeq

As in Section 5.1,
one can introduce Hamiltonian connections and Hamiltonian forms on $\bpi$. 
Let 
\beq
H=p_{it}dy^i +p_idy^i_t -\cH(t,y^i,y^i_t, p_{it},p_i)dt \label{m201}
\eeq
be a Hamiltonian form (\ref{m46}) on $\bpi$. The corresponding Hamilton
equations (\ref{m41a}) -- (\ref{m41b}) read
\bea
&& d_ty^i=\frac{\dr\cH}{\dr p_{it}}, \label{m202a}\\
&& d_ty^i_t=\frac{\dr\cH}{\dr p_i}, \label{m202b}\\
&& d_tp_i=-\frac{\dr\cH}{\dr y^i_t}, \label{m202c}\\
&& d_tp_{it}=-\frac{\dr\cH}{\dr y^i}. \label{m202d}
\eea

\begin{example}\label{j1}
Given a connection $\G$ on $\Y$, we can bring (\ref{m201}) into the form
\be
H=d_t[p_i(dy^i-\G^idt)] -\wt\cH dt =
p_{it}dy^i +p_idy^i_t -d_t(p_i\G^i)dt -\wt\cH dt, 
\ee
where $d_t\G$ is the jet lift (\ref{m190}) of $\G$ onto $J^1Y\to X$.
In particular, every Hamiltonian $\cH_\Pi$ on $\Pi=V^*Y$ defines the
Hamiltonian 
\beq
\cH=d_t\cH_\Pi=\dr_t\cH_\Pi+ y^i_t\frac{\dr\cH_\Pi}{\dr y^i} +
p_{it}\frac{\dr\cH_\Pi}{\dr p_i}
\label{m204}
\eeq
on $\bpi$ (\ref{m197}). In this case, the equations (\ref{m202a}) --
(\ref{m202d}) take the form
\bea
&& d_ty^i=\frac{\dr\cH_\Pi}{\dr p_i}, \label{m205a}\\
&& d_ty^i_t=d_t\frac{\dr\cH_\Pi}{\dr p_i}, \label{m205b}\\
&& d_tp_i=-\frac{\dr\cH_\Pi}{\dr y^i}, \label{m205c}\\
&& d_tp_{it}=-d_t\frac{\dr\cH_\Pi}{\dr y^i}. \label{m205d}
\eea
It is easily observed that they are 
equivalent to the Hamilton equations (\ref{m205a}), (\ref{m205c}) for the
Hamiltonian $\cH_\Pi$ on $\Pi$.
\end{example}

Substitution of (\ref{m202a}) into (\ref{m202b}) and of (\ref{m202c}) into
(\ref{m202d}) leads to the equations
\bea
&& d_t\frac{\dr\cH}{\dr p_{it}}=\frac{\dr\cH}{\dr p_i}, \label{m206a}\\
&& d_t\frac{\dr\cH}{\dr y^i_t}=\frac{\dr\cH}{\dr y^i} \label{m206b}
\eea
which look like the Euler--Lagrange equations for the "Lagrangian" $\cH$.
Though $\cH$ is not a true Lagrangian function, one can write 
$\cH=-\cL +d_t(p_i\G^i)$, so that 
the equations (\ref{m206a}) -- (\ref{m206b})
become the Euler--Lagrange equations for the Lagrangian $\cL$ on $\bpi$.

The solutions of the Hamilton equations (\ref{m202a}) -- (\ref{m202d})
are obviously the solutions of the Euler--Lagrange equations (\ref{m206a}) --
(\ref{m206b}), but the converse is not true. 

\begin{example}\label{j2}
Let $\cH=-\cL_Y +d_t(p_i\G^i)$, where $\cL_Y$ is a Lagrangian on $J^1Y$.
In this case, the equations  (\ref{m206a}) --
(\ref{m206b}) are equivalent to the 
Euler--Lagrange equations (\ref{m206b}) for
the Lagrangian $\cL_Y$. However, their solutions fail to be solutions of the
corresponding Hamilton equations (\ref{m202a}) -- (\ref{m202d}) in general.
\end{example}

To give a unified picture of Examples \ref{j1} and \ref{j2}, let us
consider the Hamiltonian
\beq
\cH=d_t\cH_\Pi + (p_iy^i_t -\cH_\Pi) -\cL_Y, \label{m207}
\eeq
where $\cL_Y$ is a semiregular Lagrangian on the configuration space $J^1Y$
and $\cH_\Pi$ is a Hamiltonian associated with $\cL_Y$.
The corresponding Hamilton equations (\ref{m202a}) -- (\ref{m202d}) read
\bea
&& d_ty^i=\frac{\dr\cH_\Pi}{\dr p_i}, \label{m208a}\\
&& d_ty^i_t=d_t\frac{\dr\cH_\Pi}{\dr p_i} + y^i_t-\frac{\dr\cH_\Pi}{\dr p_i},
\label{m208b}\\ 
&& d_tp_i=-\frac{\dr\cH_\Pi}{\dr y^i} -p_i +\frac{\dr\cL}{\dr y^i_t},
\label{m208c}\\ 
&& d_tp_{it}=-d_t\frac{\dr\cH_\Pi}{\dr y^i} +\frac{\dr\cH_\Pi}{\dr y^i} +
\frac{\dr\cL}{\dr y^i}.
\label{m208d}
\eea
Using the relations (\ref{m191a}) and (\ref{m192}), one can show that
solutions of the Hamilton equations (\ref{m41a}) -- (\ref{m41b}) for the
Hamiltonian $\cH_\Pi$ which live on the Lagrangian constraint space are
solutions of the equations (\ref{m208a}) -- (\ref{m208d}). 

In other words the equations (\ref{m208a}) -- (\ref{m208d}) on the
constraint subspace
\be
y^i_t =\frac{\dr\cH_\Pi}{\dr p_i}, \qquad 
p_i = \frac{\dr\cL}{\dr y^i_t}
\ee
on $\bpi$ are equivalent to the Hamilton equations (\ref{m205a}) --
(\ref{m205d}).

Now let us consider the Euler--Lagrange equations (\ref{m206a}) --
(\ref{m206b}) for the Hamiltonian (\ref{m207}). They read
\bea
&& d_ty^i-\frac{\dr\cH_\Pi}{\dr p_i}=0, \label{m209a}\\
&& d_tp_i- d_t\frac{\dr\cL_Y}{\dr y^i_t}=-\frac{\dr\cH_\Pi}{\dr y^i}
-\frac{\dr\cL_Y}{\dr y^i}. \label{m209b}
\eea
In accordance with Proposition \ref{pHL} and Corollary
\ref{j3}, every solution $s$ of the Euler--Lagrange equations for the
Lagrangian
$\cL_Y$ such that the relation (\ref{m210}) holds are solutions of the
equations (\ref{m209a}) -- (\ref{m209b}).

In particular, if the Lagrangian $\cL_Y$ is hyperregular, 
the equations (\ref{m208a}) --
(\ref{m208d}) and the equations (\ref{m209a}) -- (\ref{m209b}) are equivalent
to the corresponding Hamilton equations and the Euler--Lagrange equations
for $\cL_Y$ and the associated Hamiltonian.

\begin{example}\label{j5}
Let us consider the Hamiltonian form 
\beq
H=p_{it}(dy^i-\g^idt) +p_i(dy^i_t-\g^i_tdt) \label{m215}
\eeq
on $\bpi$, where $\g$ is the connection (\ref{m213}) on $J^1Y\to X$. 
The associated Hamilton equations (\ref{m202a}), (\ref{m202b}) read
\ben
&& d_ty^i=\g^i=y^i_t= V^*\g^i, \qquad d_ty^i_t=\g^i_t=  A^i_t
+A^i_{tj}y^j_t=V^*\g^i_t,
\label{m217}\\
&& d_tp_i =-p_j\frac{\dr\g^j_t}{\dr y^i_t} -p_{it}=V^*\g_i, \qquad
d_tp_{it} =- p_j\frac{\dr\g^j_t}{\dr y^i}=V^*\g_it, \nonumber
\een
where $V^*\g$ is the covertical connection (\ref{44}) on $\bpi=V^*J^1Y$.
The equations (\ref{m217}) recover the evolution equation (\ref{m174})
which consists with the Euler--Lagrange equation (\ref{m206a}). 
\end{example}

Turn now to the Poisson structure generated on $\bpi$  by the canonical form
$\bom$ (\ref{m200}). The corresponding Poisson bracket (\ref{m72}) reads
\beq
\{f,g\}_V =\frac{\dr f}{\dr p_it}\frac{\dr g}{\dr
y^i} + \frac{\dr f}{\dr p_i}\frac{\dr g}{\dr y^i_t} -
\frac{\dr g}{\dr p_it}\frac{\dr f}{\dr
y^i} + \frac{\dr g}{\dr p_i}\frac{\dr f}{\dr y^i_t}.
\label{m216}
\eeq
In particular, if $f$ is a function on $\Pi$ and $\cH$ is the Hamiltonian
(\ref{m204}), the Hamiltonian evolution equation (\ref{m96}) consists with
that for the Hamiltonian $\cH_\Pi$. If $f$ is a function on $J^1Y$ and $\cH$ 
is the Hamiltonian (\ref{m215}), the Hamiltonian evolution equation consists
with the evolution equation 
\be
d_{Ht}f= \dr_\g\rfloor df=\dr_tf + y^i_t\dr_if +(A^i_t
+A^i_{tj}y^j_t)\dr_i^tf.
\ee

It is readily observed that the canonical form (\ref{m200}) and the Poisson
bracket (\ref{m216}) are invariant under the canonical transformations of
$\bpi=J^1\Pi$ generated by the canonical transformations of $\Pi$.

\newpage

\section{Conservation laws and integrals of motion}

In autonomous mechanics, an integral of motion, by definition, is a functions
on the phase space whose Poisson bracket with a
Hamiltonian is equal to zero. We can not extend this description to 
time-dependent mechanics since the Hamiltonian evolution equation
(\ref{m96}) is not reduced to the Poisson bracket. 

To discover integrals of motion in time-dependent mechanics, we follow the
field theory approach, where the first
variational formula of the calculus of variations can be utilized in
order to discover differential conservation laws. This formula provides the
canonical decomposition of the Lie derivative of a Lagrangian along 
vector fields corresponding to 
infinitesimal gauge transformations into two
terms. The first one contains the variational derivatives and
vanishes on shell.  The other term is the divergence
of the corresponding symmetry flow $\cT$. If a Lagrangian is
gauge-invariant, its  Lie derivative is equal to zero and the
weak conservation law $0\ap d_\la\cT^\la$ holds on shell. 

\subsection{Lagrangian conservation laws}

In field theory, differential conservation laws are derived from
the condition of Lagrangians to be invariant under 
1-parameter groups of gauge
transformations.

By a gauge transformation is meant an isomorphism
$\Phi$ of a bundle $\pi:Y\to X$ over a diffeomorphism $f$ of $X$.
Every 1-parameter group $\Phi[\al]$ of isomorphisms of 
$Y\to X$ yields the complete vector field
\beq
u=u^\la(x^\m)\dr_\la +u^i(x^\m,y^j)\dr_i \label{02}
\eeq
which is the generator of $\Phi[\al]$. It is projected onto the vector field
$\tau=u^\m\dr_\m$ on $X$ which is the generator of 
$f[\al]$.  Conversely, one can think of any
projectable vector field (\ref{02}) on a bundle
$Y$ as being the generator of a local 1-parameter gauge group.
Using the canonical lift (\ref{1.21}) of $u$ onto $J^1Y$, we have
\beq
\bL_{\ol u}L=d(u\rfloor L)+ u\rfloor dL=[\dr_\la u^\la\cL +(u^\la\dr_\la+
u^i\dr_i +(d_\la u^i -y^i_\m\dr_\la u^\m)\dr^\la_i)\cL]\om. \label{04}
\eeq

The first variational formula provides the canonical decomposition
(\ref{C30}) of the Lie derivative 
(\ref{04}) in accordance with the variational
task. It is given by the coordinate relation
\ben
&& \dr_\la u^\la\cL +[u^\la\dr_\la+
u^i\dr_i +(d_\la u^i -y^i_\m\dr_\la u^\m)\dr^\la_i]\cL\equiv \label{C300}\\
&& \qquad   (u^i-y^i_\m u^\m )(\dr_i-d_\la \dr^\la_i)\cL -
d_\la[\pi^\la_i(u^\m y^i_\m -u^i) -u^\la\cL], \nonumber
\een
where 
\beq
\cT=\cT^\la\om_\la =[\pi^\la_i(u^\m y^i_\m-u^i )- u^\la\cL]\om_\la, 
 \quad
\pi^\la_i=\dr^\la_i\cL, \quad \om_\la=\dr_\la\rfloor\om, \label{Q30}
\eeq
is the symmetry flow along the vector field $u$.

The first variational formula (\ref{C300}) on shell (\ref{006})
comes to the weak transformation law 
\ben
&&\dr_\la u^\la\cL +[u^\la\dr_\la+
u^i\dr_i +(d_\la u^i -y^i_\m\dr_\la u^\m)\dr^\la_i]\cL  \label{J4}\\
&&\qquad \ap -d_\la[\pi^\la_i(u^\m y^i_\m -u^i)-u^\la\cL]. \nonumber
\een
If the Lie derivative $\bL_{\ol u}L$
(\ref{04}) vanishes, we have the conservation law 
\be
0\ap d_\la[\pi^\la_i(u^\m y^i_\m-u^i )-u^\la\cL]. 
\ee
It is brought into the differential conservation law
\be
0 \ap \frac{d}{dx^\la}(\pi^\la_i(u^\m \dr_\m s^i-u^i) -u^\la\cL ) 
\ee
on solutions $s$ of the Euler--Lagrange equations (\ref{2.29})

Background fields break
conservation laws as follows. Let us consider the product $Y\times Y'$
of a bundle $Y$ coordinatized by $(x^\la, y^i)$ whose sections are dynamical
fields on shell (\ref{006}) and a bundle
$Y'$ coordinatized by $(x^\la, y^A)$ whose sections are background
fields which take the background values
$y^B=\f^B(x)$, $y^B_\la= \dr_\la\f^B(x)$. Let 
\beq
u=u^\la(x)\dr_\la + u^A(x^\mu,y^B)\dr_A + u^i(x^\mu,y^B, y^j)\dr_i
\label{l68} 
\eeq
be a projectable vector field on $Y\times Y'$ which is
projected also onto $Y'$ (gauge 
transformations of background fields do not depend on the dynamic ones).
Substitution of (\ref{l68}) into (\ref{C300}) leads to
the first variational formula in the presence of background fields.
The weak identity
\be
&&\dr_\la u^\la\cL +[u^\la\dr_\la+  u^A\dr_A +
u^i\dr_i +(d_\la u^A -y^A_\m\dr_\la
u^\m)\dr^\la_A +  (d_\la u^i -y^i_\m\dr_\la u^\m)\dr^\la_i]\cL\ap  \\
&&\quad -d_\la[\pi^\la_i(u^\m y^i_\m -u^i)
-u^\la\cL]+ (u^A-y^A_\la u^\la)\dr_A\cL +\pi^\la_Ad_\la (u^A-y^A_\mu u^\mu)
\ee
holds on shell (\ref{006}).
If a total Lagrangian is gauge-invariant, we discover the transformation law
\beq
0\ap -d_\la[\pi^\la_i(u^\m y^i_\m-u^i) -u^\la\cL]+ 
(u^A-y^A_\la u^\la)\dr_A\cL + \pi^\la_A d_\la (u^A-y^A_\mu
u^\mu) \label{l70}
\eeq
in the presence of background fields.

\begin{remark}\label{nonlag}
The transformation law (\ref{l70}) can also be applied when the dynamical
equations are not Lagrangian, but are given, e.g., by local expression
\beq
(\dr_i- d_\la\dr^\la_i)\cL +F_i(t,y^j,y^j_t)=0. \label{m237}
\eeq
In this case, the transformation law reads
\beq
F_i\ap -d_\la[\pi^\la_i(u^\m y^i_\m-u^i)].  \label{m238}
\eeq
\end{remark}

\subsection{Energy-momentum conservation laws}

The transformation law (\ref{J4}) is linear in
the vector field $u$. Hence, one can consider superposition of the 
transformation laws along different vector fields. 

Every vector field $u$ on $Y$ projected onto a vector field $\tau$ on $X$
is the sum of the lift  of $\tau$ onto $Y$ and of some
vertical vector field $\vt$ on $Y$. It follows that 
every transformation law (\ref{J4}) is the
superposition of the Noether transformation law 
\be
[\vt^i\dr_i +d_\la \vt^i \dr^\la_i]\cL 
 \ap d_\la(\pi^\la_i\vt^i) 
\ee
for the Noether flow  $\cT^\la = -\pi^\la_i\vt^i$ and of the
stress-energy-momentum (SEM) transformation law \cite{cam1,giac96,sard97}.

A vector field $\tau$ on $X$ can be lifted to $Y$ only by means of a
connection on $Y$.

Let $\tau=\tau^\m\dr_\m$ be a vector field on $X$ and
$\tau_\G=\tau^\m (\dr_\m+\G^i_\m\dr_i)$
its horizontal lift onto $Y$ by a connection $\G$. The weak identity
(\ref{J4}) along $\tau_\G$ reads
\ben
&&\dr_\m\tau^\m\cL +[\tau^\m\dr_\m
+\tau^\m\G^i_\m\dr_i +(d_\la(\tau^\m\G^i_\m)
 -y^i_\m\dr_\la\tau^\m)\dr^\la_i]\cL \ap  \label{504}\\
&& \qquad - d_\la
[\pi^\la_i\tau^\m( y^i_\m- \G^i_\m)-\dl^\la_\m\tau^\m\cL], \nonumber
\een
where 
\be
\cT_\G{}^\la{}_\m  =\pi^\la_i(y_\m^i -\G^i_\m)-\dl^\la_\m\cL 
\ee
is the SEM tensor relative to the connection $\G$. 

One may choose different connections
$\G$  in order to discover SEM conservation laws.  
The SEM flows relative to $\G$
and $\G'$ differ from each other in the Noether flow along the
vertical vector field $\vt=\tau^\m(\G^i_\m-\G'^i_\m)\dr_i$. 

If the transformation law (\ref{504}) holds for any vector field $\tau$ on
$X$, we come to the system of weak equalities
\be
(\dr_\m+\G^i_\m\dr_i +d_\la\G^i_\m \dr^\la_i)\cL\ap 
-d_\la\cT_\G{}^\la{}_\m. 
\ee
For instance, if we choose the locally trivial connection $\G_0{}^i_\m=0$,
then the identity (\ref{504}) recovers the well-known transformation law 
\beq
\frac{\dr\cL}{\dr x^\m} +\frac{d}{dx^\la} \cT_0{}^\la{}_\m (s)\ap 0, \qquad
\cT_0{}^\la{}_\m (s)= \pi^\la_i\dr_\m s^i -\dl^\la_\m\cL, \label{Q40}
\eeq
of the canonical energy-momentum tensor $\cT_0$. Though it is not a 
true tensor, the
transformation law (\ref{Q40}) on solutions $s$ of differential
Euler--Lagrange equations is well-defined. 

\subsection{Hamiltonian conservation laws}

To discover the conservation laws in the framework of the Hamiltonian
formalism, we go back to Remark \ref{joint} \cite{sard97}.

Given a Hamiltonian form $H$ (\ref{4.7}) on the Legendre bundle $\Pi$, let us
consider the Lagrangian (\ref{Q3}) on $J^1\Pi$. 
One can apply the first variational formula (\ref{C300}) to this Lagrangian
in order to get the differential conservation laws in the
framework of the polysymplectic Hamiltonian formalism. 

Every projectable vector field $u$ on $Y\to X$ can be lifted to $\Pi$ as
follows:
\beq
\wt u=u^\m\dr_\m + u^i\dr_i +( - \dr_i u^j p^\la_j -\dr_\m u^\m p^\la_i
+\dr_\m u^\la p^\m_i)\dr^i_\la. \label{Q4}
\eeq
In case of the vector field $\wt u$ (\ref{Q4}) and the Lagrangian  
$L_H$ (\ref{Q3}), the first variational formula (\ref{C300}) 
on shell (\ref{3.11a}) -- (\ref{3.11b}) takes the form
\ben
&& p^\la_iy^i_\la\dr_\m u^\m
-\dr_\la(u^\la\cH) -u^i\dr_i\cH +(d_\la u^i -\dr^i_\m\cH\dr_\la u^\m)p^\la_i
\nonumber \\
&& \qquad \ap d_\la[p^\la_i(u^i-\dr^i_\m\cH u^\m) + u^\la (p^\m_i\dr^i_\m\cH
-\cH)]. \label{Q6} 
\een
If $\bL_{\wt u}L_H=0$,  
then we get the weak conservation law
\begin{equation}
0\ap =
 d_\la[p^\la_i(u^i-\dr^i_\m\cH u^\m) + u^\la (p^\m_i \dr^i_\m\cH -\cH)]\om.
\label{Q7} 
\end{equation}
On solutions $r$ of the Hamilton equations, the
weak equality (\ref{Q7}) comes to the weak differential conservation law
\be
0 \ap - \frac{d}{dx^\la}\wt\cT^\la (r)\om 
\ee
of the flow
\be
 \wt\cT^\la (r)  =
-[r^\la_i(u^i-\dr_\m^i\cH u^\m) + u^\la (r^\m_i \dr_\m^i\cH -\cH)].
\ee

The following assertion describes the relations between
differential conservation laws in Lagrangian and Hamiltonian formalisms.
 
\begin{proposition}\label{hamconslaw} Let a 
Hamiltonian form $H$ be associated
with a semiregular Lagrangian  $L$. Let $r$ be a solution
of the Hamilton equations for $H$ which lives on the Lagrangian
constraint space $Q$ and $s$ the associated solution  of the 
Euler--Lagrange equations for $L$ so that they satisfy the conditions
(\ref{2.37}). In virtue of the relations 
(\ref{2.30b}) and (\ref{Q2}), we have 
\beq
\wt\cT (r)=\cT( \wh H\circ r),\qquad
\wt\cT (\wh L\circ J^1s) =\cT(s), \label{Q10}
\eeq
where $\cT$ is the flow (\ref{Q30}).
\end{proposition}

In particular, let $\tau$ be a vector 
field on $X$ and $\tau_\G$ its horizontal
lift onto $Y\to X$ by a connection $\G$ on $Y$. We have the corresponding
flow 
\beq
\wt\cT_\G{}^\la{}_\m  =p^\la_i\dr_\m^i\wt\cH_\G 
-\dl^\la_\m(p^\nu_i\dr^i_\nu\wt\cH_\G -\wt\cH_\G), \label{Q12}
\eeq
where $\wt\cH_\G$ is the Hamiltonian  in the splitting (\ref{4.7})
of $H$ with respect to the connection $\G$. The relations
(\ref{Q10}) shows that, on the Lagrangian constraint space $Q$, the flow 
(\ref{Q12}) can be treated as the 
Hamiltonian SEM flow relative the connection $\G$.

The weak transformation law (\ref{Q6}) of the Hamiltonian SEM flow
(\ref{Q12}) takes the form
\be
-(\dr_\m +\G^j_\m\dr_j - p^\la_i\dr_j\G^i_\m\dr^j_\la)\wt\cH_\G +
p^\la_iR^i_{\la\m} \ap - d_\la\wt\cT_\G{}^\la{}_\m. 
\ee

Let us now consider the transformation law (\ref{Q6}) when the vector
field $\wt u$ on $\Pi$ is the horizontal lift 
of a vector field $\tau$ on $X$
by means of Hamiltonian connection on $\Pi\to X$ which is associated with the
Hamiltonian form $H$. We have 
\be
\wt u =\tau^\m(\dr_\mu + \dr^i_\m\cH\dr_i +\g^\la_{i\m}\dr^i_\la).
\ee
In this case, the corresponding SEM flow reads
\be
\wt\cT^\la = -\tau^\la(p^\m_i\dr^i_\m\cH -\cH), 
\ee
and the weak transformation law takes the form
\beq
-\dr_\m\cH
+d_\la(\dr^i_\m\cH p^\la_i)\ap \dr_\m(p^\la_i\dr_\la^i\cH-\cH).\label{Q61}
 \eeq
A glance at the expression (\ref{Q61}) shows that the SEM flow is not
conserved, but we can write the transformation law
\be
-\dr_\m\cH +d_\la[
 \dr^i_\m\cH p^\la_i - \dl^\la_\m(p^\nu_i\dr_\nu^i\cH-\cH)]\ap 0.
\ee
This is exactly the Hamiltonian form of the canonical energy-momentum
transformation law (\ref{Q40}) in the Lagrangian formalism.

\subsection{Integrals of motion in time-dependent mechanics}

In Lagrangian mechanics when $X$ is a 1-dimensional manifold, 
we consider conservation law along a vector field 
\beq
u=u^t\dr_t +u^i\dr_i, \qquad u^t=0,1, \label{m223}
\eeq
on $\Y$. Its jet lift (\ref{1.21}) onto $J^1Y$ reads
\be
\ol u=u^t\dr_t +u^i\dr_i +d_tu^i\dr_i^t.
\ee
In this case, the
first variational formula (\ref{C300}) takes the form
\beq
\ol u\rfloor d\cL\equiv 
(u^i-u^ty^i_t)(\dr_i-d_t \dr^t_i)\cL -d_t\cT, \label{m218}
\eeq
where 
\beq
\cT=\pi_i(u^ty^i_t -u^i) -u^t\cL \label{m225}
\eeq
is the flow along the vector field $u$. 

The first variational formula (\ref{m218}) on shell (\ref{006})
comes to the weak transformation law 
\beq
 \ol u\rfloor d\cL\ap -d_t\cT. \label{m219}
\eeq
If the Lie derivative 
\be
\bL_{\ol u}L=(ol u\rfloor d\cL) dt=
(u^t\dr_t +u^i\dr_i +d_tu^i\dr_i^t)\cL dt
\ee
vanishes, we have the conservation law 
\be
0\ap d_t[\pi_i(u^ty^i_t-u^i )-u^t\cL]. 
\ee
It is brought into the differential conservation law
\be
0 \ap \frac{d}{dt}(\pi_i(u^t \dr_t s^i-u^i) -u^t\cL ) 
\ee
on solutions $s$ of the Euler--Lagrange equations. A glance at
this expression shows that, in mechanics, the conserved flow (\ref{m225})
plays the role of a (first) integral of motion.

Every transformation law (\ref{m219}) along a
vector field $u$ (\ref{m223}) on $Y$ 
can be represented as superposition 
of the Noether transformation law along a
 vertical vector field $u$, where $u^t=0$ and  of the energy transformation
law along a horizontal lift 
\beq
\tau_\G=\dr_t+\G^i\dr_i \label{m226}
\eeq
 of the standard vector field $\dr_t$ on
$X$ by means some connection $\G$ on $\Y$ \cite{cam1,eche95}.

If $u$ is a vertical vector field, the transformation law (\ref{m219}) reads
\be
(u^i\dr_i +d_tu^i \dr^t_i)\cL \ap d_t(\pi_iu^i). 
\ee
If the Lie derivative of $L$ along $u$ is equal to zero, we have the integral
of motion $\cT=\pi_iu^i$.

\begin{example}\label{inmo1}
Let a Lagrangian $\cL$ does not depend on some coordinate $y^1$. Then, its
Lie derivative along the vertical vector field $u=\dr_1$ is equal to zero,
and the corresponding integral of motion is the momentum $\cT=\dr_1^t\cL$.
\end{example}

The transformation law (\ref{m219}) along the horizontal lift $\tau_\G$
(\ref{m226}) takes the form
\beq
(\dr_t +\G^i\dr_i +d_t\G^i\dr_i^t)\cL =- d_t(\pi_i(y^i_t -\G^i) -\cL),
\label{m227}
\eeq
where 
\beq
E_\G= \pi_i(y^i_t -\G^i) -\cL \label{m228}
\eeq
is the energy density. Obviously, 
it depends on the choice of the connection $\G$. From the physical 
viewpoint, one can treat $\G$ as a (local) reference frame,
$\dot y^i_\G= y^i_t -\G^i$ and $E_\G$ as 
the relative velocities and the energy
density respectively with regard to this reference frame.

\begin{example}\label{inmo2}
Let us put $\G=0$ that corresponds to choice of the (local) reference frame
given by the coordinates $y^i$. In this case, the energy transformation law
(\ref{m227}) takes the familiar form
\beq
\dr_t\cL =- d_t(\pi_iy^i_t -\cL),
\label{m229}
\eeq
where $y^i_t$ can be treated as velocities with respect to the above-mentined
reference frame.
\end{example}

\begin{example}\label{inmo3}
Let us consider the bundle $\R\times\R\to \R$ coordinatized by
$(t,y)$. It describes 1-dimensional motion. 
\beq
\G^i= a^it \label{C90}
\eeq
be a connection on this bundle which defines an accelerated reference frame
 with respect to the reference frame $y$. Consider the Lagrangian 
\be
L = \frac12(y^i_t-a^it)^2dt
\ee
which describes the free particle relative to the reference frame
$\G$. It is easy to see that  the energy density (\ref{m228})
relative to connection (\ref{C90}) is conserved. It is exactly  the
energy of the free particle with respect to the reference frame $\G$.
\end{example}

We now turn to the Hamiltonian mechanics.

Given a vector field (\ref{m223}), let
\beq
\wt u=u^t\dr_\m + u^i\dr_i - \dr_i u^j p_j\dr^i, \qquad u^t=0,1,\label{m234}
\eeq
be its lift onto the phase space $V^*Y$. We consider conservation
laws in time-dependent Hamiltonian mechanics along the vector fields
(\ref{m234}). As a particular case of the transformation law (\ref{Q6}), we
have
\beq
-u^t\dr_t\cH -u^i\dr_i\cH +p_id_t u^i \ap d_t(p_iu^i -u^t\cH). \label{m230} 
\eeq
 
In case of a vertical vector field $u$, this transformation law comes to the
weak equality
\be
 -u^i\dr_i\cH \ap u^id_t p_i. 
\ee
In particular, if a Hamiltonian $\cH$ is locally independent on the
coordinate $y^i$, the momentum $p_i$ is the (local) integral of motion.

The transformation law (\ref{m230}) along the horizontal lift $\tau_\G$
(\ref{m226}) takes the form
\be
-\dr_t\cH -\G^i\dr_i\cH +p_id_t \G^i \ap -d_t\wt\cH_\G. 
\ee
It follows that, in accordance with Proposition \ref{hamconslaw}, the
Hamiltonian partner of the Lagrangian energy density $E_\G$ (\ref{m228}) is
the Hamiltonian function $\wt\cH_\G$ from the splitting (\ref{m46}).
Therefore, we can treat it as the energy function with respect to the
reference frame $\G$. In particular, if
$\G^i= 0$, we get the well-known energy transformation law
\be
\dr_t\cH =d_t\cH 
\ee
which is the Hamiltonian variant of the Lagrangian law (\ref{m229}).

\end{document}